# Möbius structures, hyperbolic ends and $k$-surfaces in hyperbolic space.

Graham Smith*

**Abstract:** We present a basic introduction to the theories of Möbius structures and hyperbolic ends and we study their applications to the theory of $k$-surfaces in 3-dimensional hyperbolic space.

## 1 - Overview.

**1.1 - Hyperbolic ends and Möbius structures.** In the words of Thurston, within the family of all three-dimensional manifolds, hyperbolic three-manifolds make up "by far the most interesting, the most complex, and the most useful" class (see [28]). In this paper, we will only be concerned with two- and three-dimensional manifolds, which we will henceforth refer to simply as **surfaces** and **manifolds** respectively. In addition, in order to prevent an avalanche of unwieldy expressions, we will call a hyperbolic manifold **geometrically finite** whenever it is complete, oriented, of finite topological type and without cusps. Our aim is to present two of the main constructs used in the study of such manifolds, namely hyperbolic ends and Möbius structures.

Hyperbolic manifolds are locally modelled on three-dimensional hyperbolic space $\mathbb{H}^3$. For ease of visualisation, it is helpful to identify this space with the open unit ball $\mathbb{B}^3_1$ in $\mathbb{R}^3$ furnished with the **Beltrami-Klein metric**

$$g^{\text{BK}}_{ij} := \frac{\delta_{ij}}{(1 - \|x\|^2)} + \frac{x^i x^j}{(1 - \|x\|^2)^2}. \tag{1.1}$$

This is called the **Beltrami-Klein model** of $\mathbb{H}^3$ (see [4]). Its most useful property for our purposes is that its metric is affine equivalent to the standard euclidean metric in the sense that the geodesics of the one coincide, as sets, with the geodesics of the other. In particular, a subset $K$ of the unit ball is convex as a subset of $\mathbb{H}^3$ if and only if it is convex as a subset of $\mathbb{R}^3$.

Let $\partial_\infty \mathbb{H}^3$ denote the ideal boundary of $\mathbb{H}^3$ which, we recall, is defined to be the space of equivalence classes of complete geodesic rays in $\mathbb{H}^3$, where two such rays are deemed equivalent whenever they are asymptotic to one another (see [2]). In the Beltrami-Klein model, equivalence classes are uniquely defined by their end points, so that $\partial_\infty \mathbb{H}^3$ identifies topologically with the unit sphere $\mathbb{S}^2_1$, and the union $\mathbb{H}^3 \cup \partial_\infty \mathbb{H}^3$ likewise identifies topologically with the closed unit ball $\overline{\mathbb{B}}^3_1$.

Let $\text{PSO}_0(3,1)$ denote the group of orientation preserving isometries of $\mathbb{H}^3$. Recall that its action extends uniquely to a continuous action on the $\mathbb{H}^3 \cup \partial_\infty \mathbb{H}^3$.

**Definition 1.1.1**

*Let $S$ be a compact, oriented surface of genus at least 2, let $\Pi$ denote its fundamental group and let $\theta : \Pi \to PSO_0(3,1)$ be an injective homomorphism with discrete image. We say that $\theta$ is a **quasi-Fuchsian representation** whenever it preserves a Jordan curve in $\partial_\infty \mathbb{H}^3$. We say that a hyperbolic manifold $X$ is **quasi-Fuchsian** whenever it is isometric to the quotient of $\mathbb{H}^3$ by the image of some quasi-Fuchsian representation.*

**Remark 1.1.1.** The quasi-Fuchsian manifold $X$ is a complete hyperbolic manifold diffeomorphic to $S \times \mathbb{R}$ (see [4] and [29]).

**Remark 1.1.2.** The Jordan curve $\Gamma$ preserved by $\theta(\Pi)$ coincides with the limit set of the $\theta(\Pi)$-orbit of every point of $\mathbb{H}^3 \cup \partial_\infty \mathbb{H}^3$. In particular, $\Gamma$ is uniquely defined by this representation.

Quasi-Fuchsian manifolds are geometrically finite. In fact, they are the archetypical examples of geometrically finite hyperbolic manifolds. Of their various interesting properties, two will concern us in particular. The first is a certain natural decomposition, which is constructed as follows. Let $\theta : \Pi \to \mathrm{PSO}_0(3,1)$ be a quasi-Fuchsian representation, let $C \subseteq \partial_\infty \mathbb{H}^3$ denote the unique Jordan curve that it preserves, and let $X := \mathbb{H}^3/\theta(\Pi)$ denote the quasi-Fuchsian manifold that it defines. Let $\tilde{K}$ denote the convex hull of $C$ in $\mathbb{H}^3$ and let $\tilde{\Omega}_1$ and $\tilde{\Omega}_2$ denote the two connected components of its complement. Since $\theta(\Pi)$ preserves $\tilde{K}$, $\tilde{\Omega}_1$ and $\tilde{\Omega}_2$, their respective quotients $K$, $\Omega_1$ and $\Omega_2$ identify with subsets of $X$, and we thus obtain the decomposition

$$X := K \cup \Omega_1 \cup \Omega_2. \tag{1.2}$$

---

* Instituto de Matemática, UFRJ, Av. Athos da Silveira Ramos 149, Centro de Tecnologia - Bloco C, Cidade Universitária - Ilha do Fundão, Caixa Postal 68530, 21941-909, Rio de Janeiro, RJ - BRAZIL

Furthermore, $K$ is the minimal, closed, convex subset onto which $X$ retracts (see [4]). More generally (see [13]), every geometrically finite hyperbolic manifold decomposes in this way as the union of such a minimal, closed, convex subset, known as its **Nielsen kernel**, and finitely many unbounded open subsets, of varying topological type, known as its **ends**.

### Definition 1.1.2

*A* **height function** *over a hyperbolic manifold $Y$ is defined to be a locally strictly convex, $C^{1,1}_{\mathrm{loc}}$ function $h : Y \to ]0, \infty[$ such that*

*(1) the gradient flow lines of $h$ are unit speed geodesics; and*

*(2) for all $t > 0$, $h^{-1}([t, \infty[)$ is complete.*

*We say that a hyperbolic manifold $X$ is a* **hyperbolic end** *whenever it carries a height function.*

**Remark 1.1.3.** Height functions, whenever they exist, are unique (see Lemma 3.2.3).

**Remark 1.1.4.** We are not aware of a similar definition of hyperbolic ends having been used before in the literature. However, we will show in Chapter 3 that Definition 1.1.2 yields a rich and coherent theory, and we believe that it has the virtues over earlier definitions of being more direct and of lending itself better to potential generalisations.

Consider now the quasi-Fuchsian manifold $X$ and its three components introduced above. By standard properties of convex subsets of hyperbolic space (see [2]), the open sets $\tilde{\Omega}_1$ and $\tilde{\Omega}_2$ are both hyperbolic ends with height functions given by distance in $\mathbb{H}^3$ to $\tilde{K}$. Since the above construction is invariant under the action of $\theta(\Gamma)$, the quotients $\Omega_1$ and $\Omega_2$ are also hyperbolic ends. More generally, the connected components of the complement of the Nielsen kernel of any geometrically finite hyperbolic manifold are hyperbolic ends so that the theory of hyperbolic ends encompasses the *large scale geometry* of geometrically finite hyperbolic manifolds.

The second property of quasi-Fuchsian manifolds that interests us concerns their asymptotic geometry. Indeed, with $X$ as above, we define its **ideal boundary** $\partial_\infty X$ to be the space of equivalence classes of complete geodesic rays in $X$ which are not contained in any compact set where, again, two such rays are deemed equivalent whenever they are asymptotic to one another. The lifts of such rays are complete geodesic rays in $\mathbb{H}^3$ whose end points are not elements of $C$, so that $\partial_\infty X$ identifies with the quotient of $\partial_\infty \mathbb{H}^3 \setminus C$ under the action of $\theta(\Pi)$.

We now recall that $\partial_\infty \mathbb{H}^3$ naturally identifies with the Riemann sphere $\hat{\mathbb{C}}$ and that the action of $\mathrm{PSO}_0(3,1)$ on $\partial_\infty \mathbb{H}^3$ identifies with the action of the Möbius group $\mathrm{PSL}(2, \mathbb{C})$ on this space. This identification is immediately visible in the Beltrami-Klein model, since here the natural holomorphic structure of $\partial_\infty \mathbb{H}^3$ is none other than the structure that it inherits as a smooth, embedded submanifold of $\mathbb{R}^3$.

### Definition 1.1.3

*Let $S$ be a surface. A* **Möbius structure** *(also known as a* **flat conformal structure***) over $S$ is an atlas $\mathcal{A}$ of $S$ in $\hat{\mathbb{C}}$ all of whose transition maps are restrictions of Möbius maps. A* **Möbius surface** *is a pair $(S, \mathcal{A})$ where $S$ is a surface and $\mathcal{A}$ is a Möbius structure over $S$. In what follows, when no ambiguity arises, we will denote the Möbius surface simply by $S$.*

For each $i$, we denote $\tilde{\Sigma}_i := \partial_\infty \tilde{\Omega}_i$, so that the complement of $C$ in $\partial_\infty \mathbb{H}^3$ decomposes as

$$\partial_\infty \mathbb{H}^3 \setminus C = \tilde{\Sigma}_1 \cup \tilde{\Sigma}_2. \tag{1.3}$$

For each $i$, $\tilde{\Sigma}_i$ is trivially a Möbius surface and, since $\theta(\Pi)$ acts on $\tilde{\Sigma}_i$ by Möbius transformations, the quotient surface

$$\Sigma_i := \partial_\infty \tilde{\Omega}_i / \theta(\Pi) \tag{1.4}$$

is also a Möbius surface. In this manner, we obtain a decomposition

$$\partial_\infty X = \Sigma_1 \cup \Sigma_2 \tag{1.5}$$

of the ideal boundary of $X$ into the union of two Möbius surfaces, each homeomorphic to $S$. More generally, the ideal boundary of any geometrically finite hyperbolic manifold consists of the union of finitely many compact Möbius surfaces, one for each end, so that the theory of Möbius structures encompasses the *asymptotic geometry* of geometrically finite hyperbolic manifolds.

We underline, however, that these theories extend beyond the theory of geometrically finite hyperbolic manifolds. Indeed, it is straightforward to construct hyperbolic ends and Möbius surfaces which do not arise respectively as the ends or ideal boundaries of such manifolds. In addition, in Sections 3.4 and 3.6, we show that every hyperbolic end $X$ also has a well-defined ideal boundary, denoted by $\partial_\infty X$, given by the space of equivalence classes of complete geodesic rays in $X$, and that this ideal boundary naturally carries the structure of a Möbius surface. Conversely, in Sections 3.5 and 3.6, we show that, for every Möbius surface $S$, there exists a canonical hyperbolic end, which we denote by $\mathcal{H}S$, and which we call its **extension**, whose ideal boundary is canonically isomorphic to $S$.* For this reason, the theories of hyperbolic ends and Möbius structures are naturally developed in tandem. In fact, in contrast to the presentation of this introduction, we find that the theory of Möbius structures precedes that of hyperbolic ends, and for this reason it will be studied first in the following sections.

In Chapters 2 and 3, we comprehensively review the foundations of these theories and the relationships between them. We have chosen as our approach to derive our results using only classical tools of hyperbolic geometry, such as geodesics, spheres, horospheres, and so on, and the reader will notice certain similarities with aspects of the work [14] of Kulkarni. In particular, we have included what we believe are simpler proofs of existing results and useful generalisations of others.

Two main themes will be of particular interest to us. The first concerns the construction and properties of certain special functions which encode global geometry in a local manner. In the case of hyperbolic ends, it will be none other than the height function defined above, whose analytic properties we will establish in some detail. In the case of Möbius surfaces, it will be a $C^{1,1}_{\text{loc}}$ section of the density bundle of the surface which we call the Kulkarni-Pinkall form. This form, first studied in [15], is naturally related to the horospherical support function of immersed surfaces in $\mathbb{H}^3$ (see [7] and [23]) and for this reason constitutes a key ingredient of the a priori estimates that we will develop for certain types of immersed surfaces in Chapter 4 and which will be discussed in the following section.

The second main theme that interests us is the construction of the operators $\partial_\infty$ and $\mathcal{H}$ mentioned above. These operators allow us to pass back and forth between the families of hyperbolic ends and Möbius surfaces. In particular, they allow us to compare the geometries of different hyperbolic ends with the same ideal boundaries. We thereby obtain the following nice result. We say that a hyperbolic end $X$ is **maximal** if it cannot be isometrically embedded in a strictly larger hyperbolic end with the same ideal boundary.

**Theorem 1.1.4, Maximality**

*For every Möbius surface $S$, the extension $\mathcal{H}S$ of $S$ is, up to isometry, the unique maximal hyperbolic end with ideal boundary $S$.*

**Remark 1.1.5.** We prove Theorem 1.1.4 in Section 3.6. In the case where $S$ is compact, this result follows from the work [22] of Scannell via the natural duality between hyperbolic ends and GHMC de Sitter spacetimes (see [9]). An independent proof of the compact case was also provided by the author in [24].

**1.2 - Infinitesimal strict convexity, quasicompleteness and the asymptotic Plateau problem.** We now review some of the powerful applications of the theories of Möbius surfaces and hyperbolic ends to the study of certain types of immersed surfaces in $\mathbb{H}^3$.

---

* The extension coincides with the hyperbolic end constructed by Kulkarni-Pinkall in Section 8 of [15], where it is called the H-hull of the Möbius surface.

**Definition 1.2.1**

*An **immersed surface** in $\mathbb{H}^3$ is a pair $(S,e)$, where $S$ is an oriented surface and $e : S \to \mathbb{H}^3$ is a smooth immersion. In what follows, we denote the immersed surface sometimes by $S$ and sometimes by $e$, depending on which is more appropriate to the context.*

We first recall some standard definitions of surface theory (c.f. [5]). Let $S$ be an immersed surface. Let $\text{U}\mathbb{H}^3$ denote the bundle of unit vectors over $\mathbb{H}^3$ and let $\pi : \text{U}\mathbb{H}^3 \to \mathbb{H}^3$ denote the canonical projection. Let $N_e : S \to \text{U}\mathbb{H}^3$ denote the positively oriented unit normal vector field over $e$. The **first**, **second** and **third fundamental forms** of $e$ are respectively the symmetric bilinear forms $\text{I}_e$, $\text{II}_e$ and $\text{III}_e$ defined over $S$ such that, for every pair $\xi$, $\nu$ of vector fields over this surface,

$$\begin{aligned}
\text{I}_e(\xi,\nu) &:= \langle De \cdot \xi, De \cdot \nu \rangle, \\
\text{II}_e(\xi,\nu) &:= \langle \nabla_\xi N_e, De \cdot \nu \rangle, \text{ and} \\
\text{III}_e(\xi,\nu) &:= \langle \nabla_\xi N_e, \nabla_\nu N_e \rangle,
\end{aligned} \tag{1.6}$$

where $\nabla$ here denotes the Levi-Civita covariant derivative of $\mathbb{H}^3$. The **shape operator** of $S$ is the field $A_e$ of endomorphisms of $TS$ defined such that

$$\text{II}_e(\cdot,\cdot) =: \text{I}_e(A_e \cdot, \cdot). \tag{1.7}$$

In particular, the third fundamental form of $S$ is expressed in terms of the first fundamental form and the shape operator by

$$\text{III}_e(\cdot,\cdot) = \text{I}_e(A_e^2 \cdot, \cdot). \tag{1.8}$$

Finally, the **extrinsic** curvature of $S$ is defined by

$$K_e := \text{Det}(A_e). \tag{1.9}$$

We now restrict attention to a class of immersed surfaces to which the theories of Möbius surfaces and hyperbolic ends naturally apply. We say that an immersed surface $S$ is **infinitesimally strictly convex** (ISC) whenever its second fundamental form is everywhere positive-definite and we say that it is **quasicomplete** whenever the metric $\text{I}_e + \text{III}_e$ is complete. When these conditions hold, we associate a natural hyperbolic end and Möbius surface to $S$ as follows. First, we denote $\mathcal{E}S := S \times [0, \infty[$ and we define the function $\mathcal{E}e : \mathcal{E}S \to \mathbb{H}^3$ by

$$\mathcal{E}e(x,t) := \text{Exp}(tN_e(x)). \tag{1.10}$$

By local strict convexity of $S$, $\mathcal{E}e$ is an immersion, and we thus furnish the manifold $\mathcal{E}S$ with the unique hyperbolic structure that makes it into a local isometry. Quasicompleteness then implies that $\mathcal{E}S$ is, in fact, a hyperbolic end (see Lemma and Definition 4.1.1), which we call the **end** of $S$. In fact, $\mathcal{E}S$ is a hyperbolic end if and only if $S$ is infinitesimally strictly convex and quasicomplete.

In order to describe the natural Möbius structure associated to $S$, we require the concept of developing maps. Let $S$ be a surface and let $\phi : S \to \hat{\mathbb{C}}$ be a local diffeomorphism. For every point $x$ of $S$, there exists a neighbourhood $U$ of $x$ over which $\phi$ restricts to a diffeomorphism onto its image $V$. The set $\mathcal{A} := (U_\alpha, V_\alpha, \phi)_{\alpha \in \mathcal{A}}$ forms an atlas of $S$ in $\hat{\mathbb{C}}$ whose transition maps are trivial, and thus a fortiori Möbius. We call $\mathcal{A}$ the **pull-back** Möbius structure of $\phi$. Given any Möbius surface $S$, we say that a local diffeomorphism $\phi : S \to \hat{\mathbb{C}}$ is a **developing map** of $S$ whenever its pull-back Möbius structure is compatible with the initial Möbius structure of the surface. Not all Möbius surfaces have developing maps, and we say that a Möbius surface is **developable** whenever a developing map exists. Finally, we define a **developed** Möbius surface to be a pair $(S, \phi)$ where $S$ is a surface and $\phi : S \to \hat{\mathbb{C}}$ is a local diffeomorphism. Naturally, in this case, we furnish $S$ with the pull-back Möbius structure of $\phi$.

Now let $S$ be a quasicomplete ISC immersed surface. We define the **horizon map**[*] $\text{Hor} : \text{U}\mathbb{H}^3 \to \partial_\infty \mathbb{H}^3$ such that, for every unit speed geodesic $\gamma : \mathbb{R} \to \mathbb{H}^3$,

$$\text{Hor}(\gamma(0)) := \underset{t \to +\infty}{\text{Lim}} \gamma(t). \tag{1.11}$$

---

[*] We are not aware of any consensus in the literature concerning the name of this function. We justify our own terminology by the fact that, for any unit vector $\xi$, $\text{Hor}(\xi)$ is the ideal point towards which it points.

The **asymptotic Gauss map** of $S$ is the function $\phi_e : S \to \partial_\infty \mathbb{H}^3$ defined by

$$\phi_e := \text{Hor} \circ N_e, \tag{1.12}$$

where $N_e$ here denotes the positively oriented unit normal of $e$. It follows by infinitesimal strict convexity of $S$ that $\phi_e$ is a local diffeomorphism from $S$ into $\partial_\infty \mathbb{H}^3$ (see [2]).[†] In particular, $(S, \phi_e)$ is a developed Möbius surface which we call the **asymptotic Gauss image** of $S$.

These constructions yield a new a priori estimate for quasicomplete ISC immersed surfaces in $\mathbb{H}^3$ which we believe to be of independent interest. In order to state this estimate, we require the following parametrisation of the space of open horoballs in $\mathbb{H}^3$ by $\Lambda^2 \partial_\infty \mathbb{H}^3$. Let $y \in \partial_\infty \mathbb{H}^3$ be an ideal point. Let $B \subseteq \mathbb{H}^3$ be an open horoball centred on $y$. Let $H$ be an open half-space whose boundary is an exterior tangent to $B$ at some point. Let $D := \partial_\infty H$ denote the ideal boundary of $H$ and let $\omega(D)$ denote the area form of its Poincaré metric. It turns out that $\omega(D)(y)$ only depends on $B$. We call $\omega_y := \omega(D)(y)$ and $y$ the **asymptotic curvature** and the **asymptotic centre** of $B$ respectively, and we verify that these data define $B$ uniquely. For all $\omega_y \in \Lambda^2 \partial_\infty \mathbb{H}^3$, we henceforth denote by $B(\omega_y)$ the open horoball in $\mathbb{H}^3$ with asymptotic centre $y$ and asymptotic curvature $\omega_y$. In this manner, we obtain the desired parametrisation of the space of open horoballs in $\mathbb{H}^3$ by $\Lambda^2 \partial_\infty \mathbb{H}^3$.

Recall now that the Kulkarni-Pinkall form of a Möbius surface $S$ is a section of $\Lambda^2 S$ (see Section 2.4).

**Theorem 1.2.2, A priori estimate**

*Let $(S, e)$ be a quasicomplete ISC immersed surface in $\mathbb{H}^3$, let $\phi$ denote its asymptotic Gauss map and let $\omega$ denote the Kulkarni-Pinkall form of the developed Möbius surface $(S, \phi)$. For all $x \in S$,*

$$e(x) \in \overline{B}(\phi_* \omega(x)). \tag{1.13}$$

**Remark 1.2.1.** Theorem 1.2.2 follows immediately from Theorem 4.2.3

Finally, we apply these constructions to the study of immersed surfaces of constant extrinsic curvature. Following the work [17] of Labourie, we make the following two definitions.

**Definition 1.2.3**

*For $k > 0$, a $k$-surface is a quasicomplete, ISC immersed surface in $\mathbb{H}^3$ of constant extrinsic curvature equal to $k$. In what follows, we denote the $k$-surface sometimes by $S$ and sometimes by $e$, depending on which is more appropriate to the context.*

**Definition 1.2.4**

*Let $(S, \phi)$ be a developed Möbius surface. For $k > 0$, we say that a $k$-surface $e : S \to \mathbb{H}^3$ is a **solution** to the asymptotic Plateau problem $(S, \phi)$ whenever its asymptotic Gauss image is equal to this Möbius surface.*

In other words, Labourie's asymptotic Plateau problem concerns the unique prescription of $k$-surfaces in terms of their asymptotic Gauss images. In [17], Labourie proved various existence and uniqueness results for solutions of this problem in a more general setting than that studied here. Further existence and continuity results were also obtained by the author in [25]. There is since grown, scattered across the literature, a rich theory around the asymptotic Plateau problem, which we will review in [27].

In Chapter 4, we apply the theories of Möbius surfaces and hyperbolic ends to the study of the asymptotic Plateau problem. In particular, using the a priori estimate of Theorem 1.2.2, we obtain the following new compactness result. First, we say that two developed Möbius surfaces $(S, \phi)$ and $(S', \phi')$ are **equivalent** whenever there exists a diffeomorphism $\alpha : S \to S'$ and a Möbius map $\beta \in \text{PSL}(2, \mathbb{C})$ such that

$$\phi' \circ \alpha = \beta \circ \phi. \tag{1.14}$$

---

[†] In fact, it is not necessary for the immersed surface to be infinitesimally strictly convex for its asymptotic Gauss map to be a local diffeomorphism. It is instead sufficient that both of its principal curvatures be different to $-1$. This properties of surfaces which satisfy this condition are studied in [7] and [23].

**Theorem 1.2.5, Monotone convergence**

*Let $(S, \phi)$ be a developed Möbius surface with universal cover not equivalent to $(\hat{\mathbb{C}}, z)$, $(\mathbb{C}, z)$ or $(\mathbb{C}, Exp(z))$. Let $(\Omega_m)_{m \in \mathbb{N}}$ be a nested sequence of open subsets of $S$ which exhausts $S$. If, for $k \in ]0, 1[$ and for all $m$, $e_m : \Omega_m \to \mathbb{H}^3$ is a $k$-surface solving the asymptotic Plateau problem $(\Omega_m, \phi|_{\Omega_m})$, then $(e_m)_{m \in \mathbb{N}}$ subconverges in the $C^{\infty}_{\mathrm{loc}}$ sense over $S$ to a $k$-surface $e_{\infty} : S \to \mathbb{H}^3$ solving the asymptotic Plateau problem $(S, \phi)$.*

**Remark 1.2.2.** Theorem 1.2.5 is proven in Theorem 4.4.3.

Upon combining Theorem 1.2.5 with the existence results proven by Labourie in [17] (see also [27]), we obtain our second new result, which completely solves the asymptotic Plateau problem for $k$-surfaces in 3-dimensional hyperbolic space.

**Theorem 1.2.6, Existence and uniqueness**

*For all $k \in ]0, 1[$, and for every developed Möbius surface $(S, \phi)$ with universal cover not equivalent to $(\hat{\mathbb{C}}, z)$, $(\mathbb{C}, z)$ or $(\mathbb{C}, Exp(z))$, there exists a unique $k$-surface $e : S \to \mathbb{H}^3$ solving the asymptotic Plateau problem $(S, \phi)$.*

**Remark 1.2.3.** Theorem 1.2.6 is proven in Theorem 4.5.5.

**Remark 1.2.4.** It is an interesting open problem to determine under what conditions a $k$-surface is complete. In Appendix A, we describe an example of a non-complete $k$-surface.

**1.3 - Schwarzian derivatives.** We henceforth parametrise $k$-surfaces in such a manner that the asymptotic Gauss map is holomorphic. By reformulating the results of the preceding section in terms of the Schwarzian derivative, we then obtain nice parametrisations of the spaces of simply connected $k$-surfaces in $\mathbb{H}^3$. Indeed, let $S$ denote either the Poincaré disk $\mathbb{D}$ or the complex plane $\mathbb{C}$. A function $\phi : S \to \hat{\mathbb{C}}$ is said to be **locally conformal** whenever it is a holomorphic local diffeomorphism. The **Schwarzian derivative** (see [18]) of any such function $\phi : S \to \hat{\mathbb{C}}$ is defined by

$$D^{\mathrm{Sch}}\phi := \left(\frac{\phi''}{\phi'}\right)' - \frac{1}{2}\left(\frac{\phi''}{\phi'}\right)^2. \tag{1.15}$$

A key property of the Schwarzian derivative is that, for any locally conformal function $\phi : S \to \hat{\mathbb{C}}$ and for any Möbius map $\alpha$,

$$D^{\mathrm{Sch}}(\alpha \circ \phi) = D^{\mathrm{Sch}}\phi. \tag{1.16}$$

Furthermore, for any holomorphic function $f : S \to \mathbb{C}$, there exists a locally conformal function $\phi : S \to \hat{\mathbb{C}}$, unique up to post-composition by a Möbius map, such that

$$D^{\mathrm{Sch}}\phi = f. \tag{1.17}$$

Let $\mathrm{Hol}(S)$ denote the space of holomorphic functions over $S$ furnished with the $C^0_{\mathrm{loc}}$ topology. For all $k > 0$, let $\widetilde{\mathrm{Imm}}_k(S)$ denote the space of $k$-surfaces $e : S \to \mathbb{H}^3$ furnished with the $C^{\infty}_{\mathrm{loc}}$ topology and let $\mathrm{Imm}_k(S)$ denote its quotient under the action of post-composition by elements of $\mathrm{PSO}_0(3, 1)$. Let $\tilde{\Sigma} : \widetilde{\mathrm{Imm}}_k(S) \to \mathrm{Hol}(S)$ be the function defined such that, for every $k$-surface $e : S \to \mathbb{H}^3$,

$$\tilde{\Sigma}[e] := D^{\mathrm{Sch}}\phi_e, \tag{1.18}$$

where $\phi_e$ denotes the asymptotic Gauss map of $e$. For any $k$-surface $e \in \widetilde{\mathrm{Imm}}_k(S)$ and for any Möbius map $\alpha$,

$$\tilde{\Sigma}_{\infty}[\alpha \circ e] = D^{\mathrm{Sch}}\phi_{\alpha \circ e} = D^{\mathrm{Sch}}(\alpha \circ \phi_e) = D^{\mathrm{Sch}}\phi_e = \tilde{\Sigma}_{\infty}[e], \tag{1.19}$$

so that, for all $k$, $\tilde{\Sigma}$ descends to a continuous functional $\Sigma : \mathrm{Imm}_k(S) \to \mathrm{Hol}(S)$.

In [25], we prove an existence and uniqueness result for solutions of asymptotic Plateau problems of hyperbolic conformal type. In the present framework, this is reformulated as follows.

**Theorem 1.3.1, Hyperbolic asymptotic Plateau problem**

*For all $k \in ]0, 1[$ and for all $f \in Hol(\mathbb{D})$, there exists a unique element $e \in Imm_k(\mathbb{D})$ such that*

$$\Sigma[e] = f. \tag{1.20}$$

*Furthermore, $e$ depends continuously on $f$. In other words, $\Sigma$ defines a homeomorphism from $Hol(\mathbb{D})$ into $Imm_k(\mathbb{D})$.*

Theorem 1.2.6 now yields the corresponding result in the parabolic case.

**Theorem 1.3.2, Parabolic asymptotic Plateau problem**

*For all $k \in ]0, 1[$ and for all $f \in Hol(\mathbb{C}) \setminus \mathbb{C}$, there exists a unique element $e \in Imm_k(\mathbb{C})$ such that*

$$\Sigma[e] = f. \tag{1.21}$$

**Remark 1.3.1.** It is not known in the parabolic case whether the solution $e$ depends continuously on the data $f$. Furthermore, we de not expect this to be the case.

**Remark 1.3.2.** Interestingly, a complementary result holds in the limiting case where $k = 1$. Indeed, by a theorem of Volkov-Vladimirova and Sasaki (see [27]), $Imm_1(\mathbb{D})$ is empty and $Imm_1(\mathbb{C})$ consists only of horospheres and universal covers of cylinders of constant radius about complete geodesics. When $e$ is a horosphere, $\Sigma[e]$ vanishes and when $e$ is a universal cover of a cylinder, $\Sigma[e]$ is a non-zero constant. For this and other reasons, for $k \in ]0, 1[$, it makes sense to identify complete geodesics and ideal points of $\partial_\infty \mathbb{H}^3$ as degenerate solutions of the asymptotic Plateau problem for $f \in \mathbb{C} \setminus \{0\}$ and $f = 0$ respectively.

**1.4 - Category theory.** Our presentation will be structured around the framework of category theory. Since category theory is not commonly used in differential geometry, we recall here its basic definitions. A **category** consists of

(1) a family $\mathcal{A}$ of mathematical objects;

(2) for any two objects $X$ and $Y$ of $\mathcal{A}$, a set $\mathrm{Mor}(X, Y)$, which we call the **morphisms** from $X$ to $Y$; and

(3) for any three objects $X$, $Y$ and $Z$ in $\mathcal{A}$, a function

$$\circ : \mathrm{Mor}(X, Y) \times \mathrm{Mor}(Y, Z) \to \mathrm{Mor}(X, Z), \tag{1.22}$$

which we call **composition**, such that

(4) for any object $X$ of $\mathcal{A}$, there exists a unique element $e \in \mathrm{Mor}(X, X)$, which we call the **identity**, such that for any other object $Y$ of $\mathcal{A}$, and for all $f \in \mathrm{Mor}(X, Y)$,

$$f \circ e = e \circ f = f; \text{ and} \tag{1.23}$$

(5) for any four objects $X$, $Y$, $Z$ and $W$ of $\mathcal{A}$, for all $\alpha \in \mathrm{Mor}(X, Y)$, $\beta \in \mathrm{Mor}(Y, Z)$ and $\gamma \in \mathrm{Mor}(Z, W)$,

$$\alpha \circ (\beta \circ \gamma) = (\alpha \circ \beta) \circ \gamma. \tag{1.24}.$$

It is crucial for a correct understanding of category theory to pay close attention to the semantics of these definitions. A family is *not* a set. In fact, there is an implicit abuse of language in the concept of family: a family is a list of axioms which can be written down. Likewise, an object of a family is *not* an element of a set: it is a mathematical object which satisfies the axioms of the family. Thus, the family of groups is given by the axioms of group theory; the family of vector spaces is given by the axioms of linear algebra; and so on.

Most familiar mathematical constructs lie within this framework. For example, the category of vector spaces is the category whose objects are vector spaces and whose morphisms are linear maps; the category of Banach spaces is a category whose morphisms are bounded linear maps; the category of smooth manifolds is a category whose morphisms are smooth maps; and so on. It should hopefully become clear that in defining new mathematical objects, it is indeed often desirable to identify their morphisms and to verify whether these morphisms include identity elements and compose associatively. It is in this sense that the above axioms constitute a check-list of properties that families of mathematical objects ought to possess.

A **(covariant) functor** $\mathcal{F}$ between two categories $\mathcal{A}$ and $\mathcal{B}$ consists of

(1) a mathematical operation that associates to every object $X$ of $\mathcal{A}$ an object $\mathcal{F}(X)$ of $\mathcal{B}$; and

(2) another mathematical operation which associates to every pair $X$ and $Y$ of objects of $\mathcal{A}$ and to every morphism $\alpha$ in $\mathrm{Mor}(X, Y)$ a morphism $\mathcal{F}(\alpha)$ in $\mathrm{Mor}(\alpha(X), \alpha(Y))$,

such that

(3) for any object $X$ of $\mathcal{A}$,

$$\mathcal{F}(e) = e; \text{ and} \tag{1.25}$$

(4) for any three objects $X$, $Y$ and $Z$ of $\mathcal{A}$, for all $\alpha \in \mathrm{Mor}(X, Y)$ and for all $\beta \in \mathrm{Mor}(Y, Z)$,

$$\mathcal{F}(\beta \circ \alpha) = \mathcal{F}(\beta) \circ \mathcal{F}(\alpha). \tag{1.26}$$

Condition (4) can also be replaced with the condition that

$$\mathcal{F}(\beta \circ \alpha) = \mathcal{F}(\alpha) \circ \mathcal{F}(\beta), \tag{1.27}$$

in which case the functor is said to be **contravariant**. However, although the simplest examples of functors are often contravariant, only covariant functors will be used in this paper.

As before, it is crucial to pay close attention to the semantics of these definitions. A functor is *not* a function: it is a list of mathematical operations which can be written down. For example, the dual operation, which associates to every vector space its dual vector space is a contravariant functor from the category of vector spaces to itself; the $C^\infty$ operation, which associates to every smooth manifold the vector space of smooth functions defined over that manifold, is a contravariant functor from the category of smooth manifolds to the category of vector spaces; and so on. Once again, it should hopefully become clear that in defining new mathematical operations between families of objects, it is often desirable to know their effects on morphisms so that the above axioms again provide a check-list of properties that such operations ought to possess.

For geometers, who are used to expressing their ideas in terms of sets and functions, this formalism is often at first unsettling. However, the concepts of category theory are, ironically, *less* abstract than those of set theory and *closer* to what we have in mind when mathematical operations are performed. To see this, recall that sets are actually abstract mathematical objects which are not necessarily constructible in any sense that we would normally understand, which is precisely what gives the mystery to such results as the Banach-Tarski paradox. Families, on the other hand, are clearly defined by fixed lists of axioms which can be written down. Likewise, functions are abstract objects of set theory which are also not necessarily constructible any sense that we would normally understand, whilst functors are fixed lists of mathematical operations which can again be written down. In fact, whenever we carry out explicit calculations, we never work with functions, but rather with the sequences of mathematical operations used to define them. Such sequences, which we regularly encounter in our day-to-day mathematical life, are, in fact, closer in kind to the functors of category theory than they are to the functions of set theory.

**1.5 - Acknowledgements.** The author is grateful to Sébastien Alvarez and François Fillastre for helpful comments on earlier drafts of this paper. Figure A.1.1 was prepared by Débora Mondaini.

## 2 - Möbius structures.

**2.1 - Möbius structures.** A **Möbius structure** (also known as a **flat conformal structure**) over a surface $S$ is an atlas $\mathcal{A}$ all of whose transition maps are restrictions of Möbius maps. A **Möbius surface** is a pair $(S, \mathcal{A})$ where $S$ is a surface and $\mathcal{A}$ is a Möbius structure over this surface. In what follows, we will denote the Möbius surface simply by $S$ whenever the atlas is clear from the context. The family of Möbius surfaces forms a category whose morphisms are those functions $\phi : X \to X'$ whose expressions with respect to every pair of coordinate charts are restrictions of Möbius maps. Naturally, we identify Möbius surfaces which are isomorphic.

Every Möbius structure trivially defines a holomorphic structure over the same surface. We call the resulting Riemann surface the **underlying Riemann surface** of the Möbius surface. The operation which associates to a Möbius surface its underlying Riemann surface is trivially a covariant functor. This distinction between Möbius surfaces and their underlying Riemann surfaces is more than a mere abstract formality, and the reader may consult, for example, [6] for an overview of the rich theory concerning the relationship between the two.

The model examples of Möbius surfaces are the open subsets of $\hat{\mathbb{C}}$ and their quotients under actions of subgroups of the Möbius group $\mathrm{PSL}(2,\mathbb{C})$. More generally, given any surface $S$, and a local diffeomorphism $\phi : S \to \hat{\mathbb{C}}$, a Möbius structure is constructed over $S$ as follows. For every point $x \in S$, there exists a neighbourhood $U$ of $x$ over which $\phi$ restricts to a diffeomorphism onto its image $V$. The set $(U_\alpha, V_\alpha, \phi)_{\alpha \in A}$ of all such charts defines an atlas of $S$ in $\hat{\mathbb{C}}$ whose transition maps are trivial, and thus a fortiori Möbius. We call this structure the **pull-back structure** of $\phi$ and we denote it by $\phi^*\hat{\mathbb{C}}$. It will often be convenient in the sequel to denote the Möbius surface defined by the developing map $\phi : S \to \hat{\mathbb{C}}$ by $(S, \phi)$.

Given a Möbius surface $S$, we say that a local diffeomorphism $\phi : S \to \hat{\mathbb{C}}$ is a **developing map** of $S$ whenever its pull-back Möbius structure is compatible with the initial Möbius structure of $S$. Any two developing maps $\phi, \phi' : S \to \hat{\mathbb{C}}$ are related to one another by

$$\phi' = \alpha \circ \phi, \tag{2.1}$$

for some Möbius map $\alpha$, so that the family of all developing maps over a given Möbius surface can be parametrised by $\mathrm{PSL}(2,\mathbb{C})$ whenever it is non-empty. We say that a Möbius surface is **developable** whenever it has a developing map. In particular, every simply connected Möbius surface has this property. In the following sections, we will mainly be concerned with developable Möbius surfaces. In particular, we will take the developing maps to be given, and we leave the reader to verify that our constructions are independent of the developing maps chosen.

Non-developable Möbius surfaces are studied as follows. Given a Möbius surface $S$ with fundamental group $\Pi$ and universal cover $\tilde{S}$, any developing map $\phi$ of $\tilde{S}$ is equivariant with respect to a unique homomorphism $\theta : \Pi \to \mathrm{PSL}(2,\mathbb{C})$ which we call its **holonomy**. Furthermore, given another developing $\phi'$ with holonomy $\theta'$, there exists a unique Möbius map $\alpha$ such that

$$\begin{aligned} \theta' &= \alpha\theta\alpha^{-1}, \text{ and} \\ \phi' &= \alpha \circ \phi. \end{aligned} \tag{2.2}$$

Although non-developable Möbius surfaces will be of little interest to us in the sequel, their study has produced a deep and fascinating literature. For example, the question of which homomorphisms arise as holonomies of Möbius surfaces is addressed thoroughly by Gallo-Kapovich-Marden in [10]. Likewise, the structure of the space of Möbius surfaces with a given fixed holonomy $\theta$ is studied by Goldman in [11]. Finally, branched Möbius structures, for which the developing map is allowed to be a ramified covering, add yet another layer of sophistication to this theory (see, for example, [3])

We conclude this section by describing a key trichotomy of the theory. We say that a connected Möbius surface is **elliptic** or **parabolic** whenever its universal cover is isomorphic to $(\hat{\mathbb{C}}, z)$ or to $(\mathbb{C}, z)$ respectively and **hyperbolic** otherwise.

### Lemma 2.1.1

*Let $S$ be a connected Möbius surface. If $S$ contains an elliptic surface, then $S$ is elliptic. If $S$ contains a parabolic surface, then $S$ is either elliptic or parabolic.*

**Proof:** Upon taking universal covers, we may suppose that $S$ is simply connected. Let $S'$ be an open subset of $S$. If $S'$ is elliptic then, being compact, it is closed so that, by connectedness, $S = S'$ is also elliptic. Suppose now that $S'$ is parabolic. Let $\phi : S \to \hat{\mathbb{C}}$ be a developing map such that $\phi(S') = \mathbb{C}$. We claim that $S'$ is also simply connected. Indeed, let $\tilde{S}'$ denote its universal cover and let $\pi : \tilde{S}' \to S'$ denote the canonical projection. Since $(\phi \circ \pi)$ is a developing map of $\tilde{S}'$, it is a diffeomorphism from $\tilde{S}'$ onto $\mathbb{C}$. It follows that $\pi$ is injective and $S'$ is therefore simply connected, as asserted. In particular, $\phi$ restricts to a diffeomorphism from $S'$ onto $\mathbb{C}$.

Suppose now that $S' \neq S$. In particular, the topological boundary $\partial S'$ of $S'$ in $S$ is non-empty. Since the restriction of $\phi$ to $S'$ is a diffeomorphism, $\phi(\partial S') = \{\infty\}$. Now let $x$ be a point of $\partial S'$. Let $\Omega$ be a connected neighbourhood of $x$ in $S$ over which $\phi$ restricts to a diffeomorphism. In particular, by injectivity, $\partial S' \cap \Omega = \{x\}$. It follows that $S' \cap (\Omega \setminus \{x\})$ is a non-trivial, open and closed subset of $\Omega \setminus \{x\}$ so that, by connectedness, $\Omega \setminus \{x\} \subseteq S'$. Since $\phi(S \setminus \Omega)$ is uniformly bounded away from $\infty$, $x$ is in fact the only element of $\partial S'$. We conclude that $\phi$ defines a diffeomorphism from $S$ onto $\hat{\mathbb{C}}$, so that $S$ is elliptic. This completes the proof. $\square$

We underline that the above trichotomy for Möbius surfaces differs from the elliptic-parabolic-hyperbolic trichotomy for Riemann surfaces. Indeed, although the underlying Riemann surface of any elliptic or parabolic Möbius surface is also respectively elliptic and parabolic, there exist many hyperbolic Möbius surfaces - such as, for example, $(\mathbb{C}^*, z)$, $(\mathbb{C}, e^z)$ and $(\mathbb{C}^*, e^z)$ - whose underlying Riemann surfaces are parabolic.

**2.2 - The Möbius disk decomposition and the join relation.** We now introduce a canonical decomposition of Möbius surfaces which will be the main tool used for their study in the sequel. Let $S$ be a developable Möbius surface with developing map $\phi$. A **Möbius disk** in $S$ is a pair $(D, \alpha)$ where $D \subseteq \hat{\mathbb{C}}$ is an open disk and $\alpha : D \to S$ satisfies

$$\phi \circ \alpha = \mathrm{Id}. \tag{2.3}$$

We call the set $(D_i, \alpha_i)_{i \in I}$ of all Möbius disks in $S$ its **Möbius disk decomposition**. Since $\phi$ is a local diffeomorphism, every point of $S$ lies in the image of some Möbius disk, so that the Möbius disk decomposition covers $S$. We define the **join relation** $\sim$ of the Möbius disk decomposition such that, for all $i, j \in I$,

$$i \sim j \;\Leftrightarrow\; \alpha_i(D_i) \cap \alpha_j(D_j) \neq \emptyset. \tag{2.4}$$

This relation is trivially reflexive and symmetric, but not transitive. Composing with $\phi$, we obtain

$$i \sim j \;\Rightarrow\; D_i \cap D_j \neq \emptyset, \tag{2.5}$$

and

$$i \sim j, \; j \sim k, \; D_i \cap D_j \cap D_k \neq \emptyset \;\Rightarrow i \sim k. \tag{2.6}$$

We call the pair $((D_i)_{i \in I}, \sim)$ the **combinatorial data** of $S$. This data is sufficient to recover $S$ uniquely up to isomorphism, as follows from the following general result.

**Theorem and Definition 2.2.1**

*Let $M$ be a smooth manifold. Let $(\Omega_i)_{i \in I}$ be a family of open subsets of $M$ and let $\sim$ be a reflexive and symmetric relation over $I$ such that*

*(1) for all $i, j \in I$, $\Omega_i \cap \Omega_j$ has at most 1 connected component;*

*(2) $i \sim j \;\Rightarrow\; \Omega_i \cap \Omega_j \neq \emptyset$; and*

*(3) $i \sim j, \; j \sim k, \; \Omega_i \cap \Omega_j \cap \Omega_k \neq \emptyset \;\Rightarrow\; i \sim k$.*

*There exists a (not necessarily second-countable) smooth manifold $N$, a smooth local diffeomorphism $\phi : N \to M$ and, for all $i$, a smooth function $\alpha_i : \Omega_i \to N$ such that,*

*(A) $(\alpha_i(\Omega_i))_{i \in I}$ covers $N$;*

*(B) $i \sim j \;\Leftrightarrow\; \alpha_i(D_i) \cap \alpha_j(D_j) \neq \emptyset$; and*

*(C) for all $i$, $\phi \circ \alpha_i = \mathrm{Id}$.*

Furthermore, the triplet $(N, \phi, (\alpha_i)_{i \in I})$ is unique in the sense that if $(N', \phi', (\alpha'_i)_{i \in I})$ is another such triplet, then there exists a unique diffeomorphism $\psi : N \to N'$ such that $\phi = \phi' \circ \psi$. We call $N$ the **join** of $((\Omega_i)_{i \in I}, \sim)$, we call $\phi$ the **canonical immersion** and we call $(\alpha_i)_{i \in I}$ the **canonical parametrisations**.

**Remark 2.2.1.** If $M$ possesses any additional structure - such as, say, a hyperbolic structure, a Möbius structure, and so on - then $N$ inherits this structure from $M$, as follows immediately from the triviality of the transition maps of the atlas constructed in the proof below.

**Remark 2.2.2.** In every case where Theorem 2.2.1 will be used in the sequel, second-countability is obtainable upon covering $N$ by a countable subfamily of $(\alpha_i(\Omega_i))_{i \in I}$. However, we will never actually require this property, for two reasons. First, second-countability is only required in manifold theory for constructions involving either Sard's Theorem or partitions of unity, neither of which are used in this paper. Second, one of the more remarkable corollaries of Riemann's uniformisation theorem is that every Riemann surface is second-countable anyway, from which this property follows in all cases of interest to us.

**Proof:** We first prove existence. Define

$$\tilde{N} := \sqcup_{i \in I} \Omega_i,$$

and define the relation $\approx$ over $\tilde{N}$ such that, for all $x_i \in \Omega_i$ and $y_j \in \Omega_j$,

$$x_i \approx y_j \iff i \sim j \text{ and } x_i = y_j.$$

It follows by (3) that $\approx$ is an equivalence relation over $\tilde{N}$. Let $N := \tilde{N}/\approx$ denote its quotient space furnished with the quotient topology and let $\alpha : \tilde{N} \to N$ denote the canonical projection. Recall now that a manifold is defined to be a second-countable, Hausdorff space furnished with an atlas. The atlas of $N$ is constructed as follows. For all $i$, we verify that $\alpha$ restricts to a diffeomorphism from $\Omega_i$ onto an open subset of $N$, and we denote

$$U_i := \alpha(\Omega_i), \; V_i := \Omega_i, \; \alpha_i := \alpha|_{V_i}, \text{ and } \phi_i := \alpha_i^{-1}.$$

The family $\mathcal{A} := (U_i, V_i, \phi_i)_{i \in I}$ forms an atlas of $N$ all of whose transition maps are trivial, and thus a fortiori smooth, as desired.

Since we are not concerned with second-countability, it only remains to show that $N$ is Hausdorff. For this, let $x_i \in \Omega_i$ and $y_j \in \Omega_j$ be such that there exists a sequence $(p_m)_{m \in \mathbb{N}}$ of points in $N$ converging simultaneously to $\alpha(x_i)$ and to $\alpha(y_j)$. For sufficiently large $m$, $p_m$ has representative elements $x_{m,i}$ in $\Omega_i$ and $y_{m,j}$ in $\Omega_j$ respectively, which converge towards $x_i$ and $y_j$ respectively. In particular, $i \sim j$ and, for all $m$, $x_{m,i} = y_{m,j}$. Upon taking limits, we obtain $x_i = y_j$, so that $x_i \approx y_j$ and therefore $\alpha(x_i) = \alpha(y_j)$. We conclude that $N$ is indeed Hausdorff, and therefore a (not necessarily second-countable) manifold.

Finally, the canonical inclusion $\tilde{\phi} : \tilde{N} \to M$ trivially descends to a local diffeomorphism $\phi : N \to M$. We verify that $(N, \phi, (\alpha_i)_{i \in I})$ has the desired properties, thus proving existence.

To prove uniqueness, let $(N', \phi', (\alpha'_i)_{i \in I})$ be another such triplet. Define $\tilde{\psi} : \tilde{N} \to N'$ such that, for all $x_i \in \Omega_i$,

$$\tilde{\psi}(x_i) := \alpha'_i(x_i).$$

We first show that $\tilde{\psi}$ descends to a function $\psi : N \to N'$. Indeed, let $x_i \in \Omega_i$ and $y_j \in \Omega_j$ be such that $x_i \approx y_j$. By $(B)$,

$$\alpha'_i(\Omega_i) \cap \alpha'_j(\Omega_j) \neq \emptyset.$$

Furthermore, by $(1)$, $(C)$ and a connectedness argument

$$\alpha'_i|_{\Omega_i \cap \Omega_j} = \alpha'_j|_{\Omega_i \cap \Omega_j}.$$

In particular, $\tilde{\psi}(x_i) = \tilde{\psi}(y_j)$ so that $\tilde{\psi}$ indeed descends to a function $\psi : N \to N'$. By $(A)$, $\psi$ is surjective, by $(B)$ and $(C)$, it is injective. Since $\alpha_i$ and $\alpha'_i$ are local diffeomorphisms for all $i$, it follows that $\psi$ is a diffeomorphism. Finally, by $(C)$ again $\phi' \circ \psi = \phi$. This proves existence of $\psi$, and since uniqueness is trivial, this completes the proof. $\square$

**2.3 - Geodesic arcs and convexity.** We now introduce a concept of geodesics for sets of Möbius disks in a given non-elliptic Möbius surface. This in turn yields a concept of convexity for such sets which will be useful for establishing uniqueness in the constructions that follow.

To begin with, we study the geometry of the space $\mathcal{D}$ of disks in $\hat{\mathbb{C}}$. Recall that $\hat{\mathbb{C}}$ naturally identifies with the ideal boundary $\partial_\infty \mathbb{H}^3$ of $\mathbb{H}^3$. With this identification, every disk $D$ in $\hat{\mathbb{C}}$ is the ideal boundary of a unique open half-space $H$ in $\mathbb{H}^3$. The boundary $\partial H$ of every open half-space in $\mathbb{H}^3$ is a totally geodesic plane which we orient so that its positively oriented normal points outwards from $H$. Trivially, open half-spaces in $\mathbb{H}^3$ are uniquely defined by their oriented boundaries. Consequently, any parametrisation of the space of oriented totally geodesic planes in $\mathbb{H}^3$ is also a parametrisation of $\mathcal{D}$.

The space of oriented totally geodesic planes in $\mathbb{H}^3$ is parametrised by $(2,1)$-dimensional de Sitter space $\mathrm{dS}^{2,1}$ as follows. First, we identify $\mathbb{H}^3$ and $\mathrm{dS}^{2,1}$ with subsets of $\mathbb{R}^{3,1}$, namely

$$
\begin{aligned}
\mathbb{H}^3 &:= \left\{ x \in \mathbb{R}^{3,1} \mid \langle x,x \rangle_{3,1} = -1,\ x_4 > 0 \right\},\ \text{and} \\
\mathrm{dS}^{2,1} &:= \left\{ x \in \mathbb{R}^{3,1} \mid \langle x,x \rangle_{3,1} = 1 \right\},
\end{aligned}
\tag{2.7}
$$

where here $\langle \cdot, \cdot \rangle_{3,1}$ denotes the **Minkowski metric** with signature $(3,1)$, that is

$$
\langle x,x \rangle_{3,1} := x_1^2 + x_2^2 + x_3^2 - x_4^2.
\tag{2.8}
$$

With this identification, every oriented totally geodesic plane $P$ in $\mathbb{H}^3$ is the intersection of $\mathbb{H}^3$ with a unique oriented, time-like, linear hyperplane $\hat{P}$ in $\mathbb{R}^{3,1}$. Every such hyperplane has, in turn, a well-defined positively-oriented unit normal vector $N$. Since $N$ is also spacelike, it is an element of $\mathrm{dS}^{2,1}$. This yields a bijection between the space of oriented, totally geodesic planes in $\mathbb{H}^3$ and $\mathrm{dS}^{2,1}$, which is the desired parametrisation.

Recall now that a subset $\Gamma$ of $\mathrm{dS}^{2,1}$ is a geodesic if and only if it is the intersection of $\mathrm{dS}^{2,1}$ with a linear plane $\hat{\Gamma}$. Furthermore, $\Gamma$ is said to be **spacelike**, **lightlike** or **timelike** respectively whenever the restriction to this plane of the Minkowski metric has signature $(2,0)$, $(1,0)$ or $(1,1)$. Of particular interest to us will be the spacelike geodesics. Observe first that any two distinct totally geodesic planes in $\mathbb{H}^3$ with non-trivial intersection meet along a complete geodesic.

**Lemma 2.3.1**

*Let $P$ and $P'$ be distinct, oriented totally-geodesic planes in $\mathbb{H}^3$ which are neither equal nor equal with opposing orientations. $P$ and $P'$ have non-trivial intersection if and only if their corresponding points in $dS^{2,1}$ lie along a common spacelike geodesic $\Gamma$. Furthermore, motion at constant speed along $\Gamma$ corresponds to rotation at constant angular speed around their common geodesic $G$.*

**Proof:** Observe first that the orthogonal complement in $\mathbb{R}^{3,1}$ of any timelike linear plane $\hat{G}$ is a spacelike linear plane $\hat{\Gamma}$ whose intersection with $\mathrm{dS}^{2,1}$ is a circle in $\hat{\Gamma}$ and a spacelike geodesic $\Gamma$ in $\mathrm{dS}^{2,1}$. Now let $P = \hat{P} \cap \mathbb{H}^3$ and $P' = \hat{P}' \cap \mathbb{H}^3$ be oriented totally-geodesic planes in $\mathbb{H}^3$ which are neither equal nor equal with opposing orientations. These planes meet along a common geodesic $G$ if and only if $\hat{P}$ and $\hat{P}'$ contain a common timelike linear plane $\hat{G}$. This in turn holds if and only if their unit normals $N$ and $N'$, which are already elements of $\mathrm{dS}^{2,1}$, are also elements of the orthogonal complement $\hat{\Gamma}$ of $\hat{G}$. This proves the first assertion. Since the second assertion is straightforward, this completes the proof. $\square$

We now return to the case of disks in $\hat{\mathbb{C}}$. We say that two distinct disks $D_0$ and $D_1$ **overlap** whenever their boundary circles meet at exactly two points. Observe that this holds if and only if their intersection is non-trivial, the intersection of their complements is non-trivial, and neither is contained within the other. With the preceding parametrisation, this is precisely the requirement for their corresponding points in $\mathrm{dS}^{2,1}$ to lie along a common spacelike geodesic. In addition, the corresponding point of a third disk $D$ lies along the *shorter* geodesic arc between these two points if and only if

$$
D_0 \cap D_1 \subseteq D \subseteq D_0 \cup D_1.
\tag{2.9}
$$

We thus define the **geodesic arc** between two overlapping disks $D_0$ and $D_1$ to be the set of all disks $D$ in $\hat{\mathbb{C}}$ which satisfy this property.

This concept of geodesic arc extends to the Möbius disk decomposition of $S$ as follows. We say that two distinct Möbius disks $(D_0, \alpha_0)$ and $(D_1, \alpha_1)$ **overlap** whenever $\alpha_0(D_0)$ and $\alpha_1(D_1)$ have non-trivial intersection and neither is contained within the other. Upon composing with $\phi$, it follows that $D_0$ and $D_1$ likewise have non-trivial intersection and neither is contained within the other. In addition, since $S$ is not elliptic, the complements of $D_0$ and $D_1$ also have non-trivial intersection, so that $D_0$ and $D_1$ also overlap. Using a connectedness argument, we show that

$$\alpha_0|_{D_0 \cap D_1} = \alpha_1|_{D_0 \cap D_1}, \tag{2.10}$$

so that these functions join to define a function $\alpha_{01} : D_0 \cup D_1 \to S$ such that

$$\phi \circ \alpha_{01} = \mathrm{Id}. \tag{2.11}$$

In particular, for any other disk $D$ along the geodesic arc from $D_0$ to $D_1$, $(D, \alpha_{01})$ is also a Möbius disk in $S$. We thus define the **geodesic arc** from $(D_0, \alpha_0)$ to $(D_1, \alpha_1)$ to be the set of all Möbius disks in $S$ of this form. We say that any subset $(D_i, \alpha_i)_{i \in J}$ of the Möbius disk decomposition of $S$ is **convex** whenever it contains the geodesic arc between any two of its overlapping disks.

**2.4 - The Kulkarni-Pinkall form.** In [15], Kulkarni-Pinkall construct for any Möbius surface of hyperbolic type a canonical metric which encodes its global geometry in a local manner. Kulkarni-Pinkall's construction will play a central role in the $C^0$ estimates that we will derive in Chapter 4 for quasicomplete ISC immersions in $\mathbb{H}^3$. However, we will adopt here a slightly different perspective, since we believe it to be more natural to work in terms of 2-forms rather than in terms of metrics.

Let $S$ be a developable Möbius surface with developing map $\phi$, let $(D_i, \alpha_i)_{i \in I}$ denote its Möbius disk decomposition and, for all $i \in I$, let $H_i$ denote the open half-space in $\mathbb{H}^3$ with ideal boundary $D_i$. For all $x \in S$, let $I(x)$ denote the set of indices $i$ such that $x \in \alpha_i(D_i)$. For any disk $D \in \hat{\mathbb{C}}$, let $\omega(D)$ denote the area form of its unique hyperbolic metric. We define $\omega_\phi$, the **Kulkarni-Pinkall form** of $S$, such that, for all $x \in S$,

$$\omega_\phi(x) := \inf_{i \in I(x)} \phi^* \omega(D_i), \tag{2.12}$$

and we define $g_\phi$ the **Kulkarni-Pinkall metric** of $S$ by

$$g_\phi := \omega_\phi(\cdot, J\cdot), \tag{2.13}$$

where $J$ here denotes the complex structure of $S$.

**Lemma 2.4.1, Monotonicity**

*Let $S$ and $S'$ be developable Möbius surfaces with respective developing maps $\phi$ and $\phi'$ and respective Kulkarni-Pinkall forms $\omega_\phi$ and $\omega_{\phi'}$. If $\alpha : S \to S'$ is a morphism such that $\phi = \phi' \circ \alpha$, then*

$$\omega_\phi \geq \alpha^* \omega_{\phi'}. \tag{2.14}$$

**Proof:** Indeed, composition with $\alpha$ sends the Möbius disk decomposition of $S$ into the Möbius disk decomposition of $S'$. $\square$

The following family of partial orders over $I$ will prove useful in deriving properties of the Kulkarni-Pinkall form. For all $x \in S$, we define

$$i \geq_x j \iff i, j \in I(x) \text{ and } \omega(D_i)(y) \leq \omega(D_j)(y), \tag{2.15}$$

where $y := \phi(x)$. The geometric significance of the Kulkarni-Pinkall form as well as this partial order becomes clear once we recall the parametrisation of the space of open horoballs in $\mathbb{H}^3$ by $\Lambda^2 \partial_\infty \mathbb{H}^3$ described in Section 1.2. Indeed, for all $x \in S$, $\phi_* \omega_\phi(x)$ is simply the infimal asymptotic curvature of horoballs centred on $\phi(x)$ and contained in $H_i$, as $i$ varies over $I(x)$. Likewise, for all $x \in S$ and for all $i, j \in I(x)$, $i \geq_x j$ if and only if every open horoball asymptotically centred on $\phi(x)$ and contained in $H_j$ is also contained in $H_i$.

**2.5 - Analytic properties of the Kulkarni-Pinkall form.** We restrict our attention initially to the simpler case of Möbius surfaces of the form $(\Omega, z)$, where $\Omega$ is an open subset of $\hat{\mathbb{C}}$. At this stage, it is useful to recall that, for a disk $D$ in the complex plane $\mathbb{C}$ of radius $R$ with centre lying at distance $r < R$ from the origin,

$$\omega(D)(0) = \frac{4R^2 dxdy}{(R-r)^2(R+r)^2}. \tag{2.16}$$

In particular, if $\omega(D)(0) < \lambda^2 dxdy$, then $D$ contains a disk of radius $1/\lambda$ about the origin.

**Lemma and Definition 2.5.1**

*Let $\Omega$ be an open subset of $\hat{\mathbb{C}}$ and let $\omega$ denote its Kulkarni-Pinkall form.*

*(1) If the complement of $\Omega$ in $\hat{\mathbb{C}}$ contains at most 1 distinct point then, for all $x$, $\omega(x) = 0$ and $I(x)$ contains no maximal element with respect to $\geq_x$.*

*(2) If the complement of $\Omega$ in $\hat{\mathbb{C}}$ contains at least 2 distinct points then, for all $x$, $\omega(x) > 0$ and $I(x)$ contains a unique maximal element with respect to $\geq_x$ which realizes $\omega(x)$.*

*In the second case, we denote by $\max(x)$ the unique maximal element of $I(x)$.*

**Proof:** The first assertion is trivial. To prove the second assertion, we may suppose that $\Omega$ is a proper subset of the complex plane $\mathbb{C}$. Existence follows by compactness of the set of (possibly ideal) disks in $\mathbb{C}$ which have radius bounded below, which contain a fixed point $z_0$ and which avoid another fixed point $w_0$. Observe now that $\omega(D_i)(x)$ restricts to a strictly concave function over every geodesic arc in $I(x)$. Uniqueness thus follows by convexity of $I(x)$. Finally, since $\omega(x)$ is realized by the unique maximal element of $I(x)$, $\omega(x) > 0$, and this completes the proof. $\square$

Given an ideal point $x \in \partial_\infty \mathbb{H}^3$ and a closed subset $Y \subseteq \mathbb{H}^3 \cup \partial_\infty \mathbb{H}^3$, we define the **curvature of distance** $c(x, Y)$ from $x$ to $Y$ to be the infimal asymptotic curvature of open horoballs with asymptotic centre $x$ which do not meet $Y$.

**Lemma 2.5.2**

*Let $\Omega$ be a proper open subset of the complex plane $\mathbb{C}$, let $\omega$ denote its Kulkarni-Pinkall form, let $(D_i, \alpha_i)_{i \in I}$ denote its Möbius disk decomposition, and, for all $i \in I$, let $H_i$ denote the open half-space in $\mathbb{H}^3$ with ideal boundary $D_i$. Let $K$ denote the convex hull in $\mathbb{H}^3 \cup \partial_\infty \mathbb{H}^3$ of the complement of $\Omega$ and let $\pi : \Omega \to \partial K$ denote the canonical projection. For all $x \in \Omega$,*

$$\omega(x) = c(x, K), \tag{2.17}$$

*and $H(x) := H_{\max(x)}$ is the unique supporting open half-space of $K$ at the point $\pi(x)$ such that $\partial H(x)$ is orthogonal to the geodesic joining $\pi(x)$ to $x$. In particular, $\omega(x)$, $H(x)$ and $D(x) := D_{\max(x)}$ are $C^{0,1}_{loc}$ functions over $\Omega$.*

**Remark 2.5.1.** In fact, Kulkarni-Pinkall show in [15] that $\omega$ is a $C^{1,1}_{loc}$ function.

**Proof:** Since the complement of $K$ in $\mathbb{H}^3$ is the union of all open half-spaces with ideal boundary in $\Omega$, we have

$$K^c = \underset{i \in I}{\cup} H_i,$$

from which it follows that $\omega(x) = c(x, K)$ for all $x \in \Omega$. Now choose $x \in \Omega$. Let $B$ be the open horoball in $\mathbb{H}^3$ with asymptotic centre $x$ and asymptotic curvature $c(x, K)$. Since $H_{\max(x)}$ is the unique open half-space in $K^c$ containing $B$, the second assertion follows and this completes the proof. $\square$

We now address the general case. Let $S$ be a developable Möbius surface with developing map $\phi$ and let $(D_i, \alpha_i)_{i \in I}$ denote its Möbius disk decomposition. For all $x \in S$, with $I(x)$ defined as in Section 2.4, we define

$$\Omega_x := \underset{i \in I(x)}{\cup} D_i. \tag{2.18}$$

For all $i, j \in I(x)$, $\alpha_i$ coincides with $\alpha_j$ over $D_i \cap D_j$ so that the join of these functions yields a function $\alpha_x : \Omega_x \to S$ satisfying $\phi \circ \alpha_x = \mathrm{Id}$. We call $(\Omega_x, \alpha_x)$ the **localisation** of $S$ at $x$. The following trichotomy follows immediately from Lemma 2.1.1.

**Lemma 2.5.3**

*Let $S$ be a developable Möbius surface with developing map $\phi : S \to \hat{\mathbb{C}}$.*

*(1) If $S$ is elliptic, then $\Omega_x = \hat{\mathbb{C}}$ for all $x$.*

*(2) If $S$ is parabolic, then $\Omega_x$ is the complement of a single point in $\hat{\mathbb{C}}$ for all $x$.*

*(3) If $S$ is hyperbolic, then the complement of $\Omega_x$ contains at least two points in $\hat{\mathbb{C}}$ for all $x$.*

For all $x \in S$, we define $\omega_{\phi,x}$, the **local Kulkarni-Pinkall form** of $S$ at $x$, to be the push-forward through $\alpha_x$ of the Kulkarni-Pinkall form of $(\Omega_x, z)$. Since composition with $\alpha_x$ sends the Möbius disk decomposition of $(\Omega_x, z)$ to $I(x)$, Lemmas 2.5.1 and 2.5.3 immediately yield the following result.

**Lemma and Definition 2.5.4**

*Let $S$ be a developable Möbius surface of hyperbolic type with developing map $\phi : S \to \hat{\mathbb{C}}$ and Kulkarni-Pinkall form $\omega_\phi$. For all $x \in S$, $I(x)$ has a unique maximal element which realises $\omega_\phi(x)$. Furthermore*

$$\omega_\phi(x) = \omega_{\phi,x}(x), \tag{2.19}$$

*and, for all $y \in \alpha_x(\Omega_x)$,*

$$\omega_\phi(y) \leq \omega_{\phi,x}(y). \tag{2.20}$$

*For all $x \in S$, we denote by $\max(x)$ the unique maximal element of $I(x)$.*

Analytic properties of $\omega_\phi$ analogous to those obtained in Lemma 2.5.2 for localised Möbius structures follow upon refining (2.19) to equality over a neighbourhood of $x$.

**Lemma 2.5.5**

*Let $S$ be a developable Möbius surface of hyperbolic type with developing map $\phi$. For all $x \in S$, there exists a neighbourhood $U_x$ of $x$ such that, for all $y \in U_x$,*

$$\max(y) \in I(x). \tag{2.21}$$

**Proof:** Since $S$ is hyperbolic, we may suppose that $\Omega_x$ is a proper subset of $\mathbb{C}$. Let $(D_i, \alpha_i)_{i \in I}$ denote the Möbius disk decomposition of $S$. Denote $i := \max(x)$. We may suppose that $D_i$ is the upper half-space in $\mathbb{C}$ and that $\phi(x) = \sqrt{-1}$. Let $d_h$ denote the hyperbolic distance in $D_i$ and define $U_x$ by

$$U_x := \left\{ y \in \alpha_i(D_i) \mid d_h(\phi(y), \phi(x)) < \mathrm{Log}((1+\sqrt{5})/2) \right\}.$$

Let $y$ be an element of $U_x$. Observe that $\phi(y)$ is contained in the euclidean ball of radius $(\sqrt{5}-1)/2$ about $\phi(x)$ in $\mathbb{C}$. Denote $j := \max(y)$. Since $S$ is hyperbolic, $\partial D_i$ intersects $\partial D_j$ at at least one point and, upon applying a suitable Möbius transformation, we may suppose that one of these points lies at infinity. In particular $D_j$ is a disk in $\mathbb{C}$. However, by (2.16) and (2.20),

$$\omega(D_j)(\phi(y)) = \phi_* \omega_\phi(y) \leq \phi_* \omega_{\phi,x}(y) \leq \frac{4dxdy}{(\sqrt{5}-1)^2}.$$

It follows by (2.16) again that $D_j$ contains a ball of radius $(\sqrt{5}-1)/2$ about $\phi(y)$. In particular, $\phi(x) \in D_j$, so that $x \in \alpha_j(D_j)$ and $j \in I(x)$, as desired. $\square$

**Corollary 2.5.6**

*Let $S$ be a developable Möbius surface of hyperbolic type with developing map $\phi : S \to \hat{\mathbb{C}}$ and Kulkarni-Pinkall form $\omega_\phi$. Let $x$ be a point of $S$. With $U_x$ as in Lemma 2.5.5, for all $y \in U_x$,*

$$\omega_\phi(y) = \omega_{x,\phi}(y). \tag{2.22}$$

Combining the above results yields a description of the analytic properties of the Kulkarni-Pinkall form of every Möbius surface.

**Theorem 2.5.7**

*Let $S$ be a developable Möbius surface with developing map $\phi : S \to \hat{\mathbb{C}}$, Kulkarni-Pinkall form $\omega_\phi$ and Möbius disk decomposition $(D_i, \alpha_i)_{i \in I}$.*

*(1) If $S$ is of elliptic or parabolic type, then $\omega_\phi$ vanishes identically.*

*(2) If $S$ is of hyperbolic type, then $\omega_\phi$ is a nowhere vanishing section of $\Lambda^2 S$. Furthermore $\omega_\phi(x)$ and $D(x) := D_{\max(x)}$ are $C^{0,1}_{\mathrm{loc}}$ functions over $S$.*

Finally, by Lemma 2.5.5 the Kulkarni-Pinkall metric of any Möbius surface of hyperbolic type is everywhere non-degenerate. In addition, the proof of this result yields the following global information concerning this metric.

**Lemma 2.5.8**

*Let $S$ be a developable Möbius surface with developing map $\phi$. The Kulkarni-Pinkall metric $g_\phi$ of $S$ is complete.*

**Proof:** It suffices to show that there exists $r_0 > 0$ such that the closed ball of radius $r_0$ with respect to $g_\phi$ about any point of $S$ is compact. Let $(D_i, \alpha_i)_{i \in I}$ denote the Möbius disk decomposition of $S$ and for all $i$, let $H_i$ denote the open half-space in $\mathbb{H}^3$ with ideal boundary $D_i$. We consider first the case where $(S, \phi) = (\Omega, z)$ for some connected neighbourhood $\Omega$ of $0$ in $\mathbb{C}$. We identify $\mathbb{H}^3$ with the upper half-space in $\mathbb{C} \times \mathbb{R}$. Let $K$ denote the convex hull in $\mathbb{H}^3$ of the complement of $\Omega$ in $\hat{\mathbb{C}}$. We may suppose that $D := D_{\max(0)}$ is the unit disk about the origin so that $(0, 1)^t$ is a boundary point of $K$. In particular, for all $j \in J$, $(0, 1)^t \notin H_j$, from which it follows that $g_\phi$ is bounded below by the standard spherical metric of $\mathbb{C}$, that is

$$g_\phi \geq g^{\mathrm{Sph}} := \frac{4}{(1 + |z|^2)^2} \delta_{ij}.$$

Consider now the general case. Let $x$ be a point of $S$. Let $(\Omega_x, \alpha_x)$ denote the localisation of $S$ at this point. As before, we may suppose that $\phi(x) = 0$ and that $D := D_{\max(x)}$ is the unit disk in $\mathbb{C}$ about the origin. Let $U_x$ be as in Lemma 2.5.5. Since the hyperbolic metric of $D$ is given by

$$g^{\mathrm{Hyp}} := \frac{4}{(1 - |z|^2)^2} \delta_{ij},$$

it follows that $\phi(U_x)$ is the euclidean disk of radius $(\sqrt{5} - 2)$. However, by the preceding paragraph, over this disk,

$$\phi_* g_\phi \geq g^{\mathrm{Sph}},$$

so that $U_x$ contains the open disk of radius $\arcsin((\sqrt{5} + 2)/10)$ with respect to $g_\phi$ about $x$. The result now follows with $r_0$ equal to this radius. $\square$

## 3 - Hyperbolic ends.

**3.1 - Hyperbolic ends.** Given a hyperbolic manifold $X$, we define a **height function** over $X$ to be a strictly convex $C^{1,1}_{\mathrm{loc}}$ function $h : X \to ]0, \infty[$ such that

(1) the gradient flow lines of $h$ are unit speed geodesics; and

(2) for all $t > 0$, $h^{-1}([t, \infty[)$ is complete.

We will see in Lemma 3.2.3, below, that height functions, whenever they exist, are unique. We define a **hyperbolic end** to be a hyperbolic manifold which carries a height function. The family of hyperbolic ends forms a category whose morphisms are those functions $\psi : X \to X'$ which are local isometries. Naturally, we identify hyperbolic ends which are isometric.

We first identify various components of hyperbolic ends. Let $X$ be a hyperbolic end with height function $h$. We call the gradient flow lines of $h$ **vertical lines**. These curves form a geodesic foliation of $X$ which we denote by $\mathcal{V}$ and which we call its **vertical line foliation**. We call the level sets of $h$ the **levels** of $X$. These form another foliation of $X$ by $C^{1,1}_{\mathrm{loc}}$ embedded surfaces which we call its **level set foliation** and which we denote by $(X_t)_{t>0}$. These two foliations are transverse to one another and every vertical line intersects every level at exactly one point. From this it follows that every level of $X$ is naturally diffeomorphic to the leaf space of $\mathcal{V}$. For all $t > 0$, we define the **vertical projection** $\pi_t : X \to X_t$ to be the function which sends each point $x$ of $X$ to the intersection with $X_t$ of the vertical line on which it lies. By standard properties of convex sets in $\mathbb{H}^3$, for all $t > 0$, $\pi_t$ restricts to a 1-Lipschitz function from $h^{-1}([t, \infty[)$ into $X_t$.

We call any local isometry $\phi : X \to \mathbb{H}^3$ a **developing map** of $X$. Any two developing maps $\phi, \phi' : X \to \mathbb{H}^3$ are related by

$$\phi' = \alpha \circ \phi, \tag{3.1}$$

for some Möbius map $\alpha$. We say that $X$ is **developable** whenever it has a developing map. In particular, every simply connected hyperbolic end has this property. In the following sections, we will only be concerned with developable hyperbolic ends and we leave the reader to formulate the trivial extensions of our results to the general case. In particular, we will take the developing maps to be given, and we leave the reader to verify that our constructions are independent of the developing maps chosen.

The model examples of hyperbolic ends are the complements $\Omega$ of closed, convex subsets $K$ of $\mathbb{H}^3$, where the height function is the distance to $K$. More sophisticated examples are given by quotients of such subsets by subgroups of the Möbius group $\mathrm{PSL}(2, \mathbb{C})$, such as the ends of quasi-Fuchsian manifolds studied in the introduction. We recall that the complement of the Nielsen kernel of every finite geometry hyperbolic manifold is the union of finitely many hyperbolic ends (see, for example, [13]). However, we emphasize again that most hyperbolic ends do not arise in this manner. Indeed, the developing map of the universal cover of any end of any finite geometry hyperbolic manifold with fundamental group not equal to $\mathbb{Z}$ is an embedding in $\mathbb{H}^3$. However, as we will see in Section 3.5, it is straightforward to construct simply connected hyperbolic ends with non-injective developing maps.

The key to understanding hyperbolic ends lies in the following analogue of the Hopf-Rinow Theorem.

**Theorem 3.1.1**

Let $(X, h)$ be a hyperbolic end. If $\gamma : [0, a[ \to X$ is a geodesic segment such that $\dot\gamma(0)$ is not downward-pointing, then $\gamma$ extends to a geodesic ray defined over the entire half-line $[0, \infty[$.

**Remark 3.1.1.** A suitably modified version of Theorem 3.1.1 holds under the weaker condition that there exists a convex function $f : X \to ]0, \infty[$ such that $f^{-1}([t, \infty[)$ is complete for all $t > 0$. In fact, using the arguments of the following sections, we may show that a hyperbolic manifold $X$ is a hyperbolic end whenever there exists a $C^{1,1}_{\mathrm{loc}}$ convex function $f : X \to ]0, \infty[$ such that $f^{-1}([t, \infty[)$ is complete for all $t > 0$ and $\|\nabla f\| \geq \epsilon > 0$. Such functions, which one may call **generalised height functions** are thus natural objects of study in the theory of hyperbolic ends (c.f. [1]).

**Proof:** By strict convexity of $h$, $(h \circ \gamma)$ has strictly increasing derivative. Since, by hypothesis, $(h \circ \gamma)$ has non-negative derivative at 0, it follows that its derivative is strictly positive for all positive time, so that $(h \circ \gamma)$ is itself strictly increasing. In particular, $\gamma$ remains within a complete subset of $X$ and may thus be extended indefinitely, as desired. $\square$

**3.2 - The half-space decomposition.** Let $X$ be a developable hyperbolic end with height function $h$ and developing map $\phi : X \to \mathbb{H}^3$. We define a **half-space** in $X$ to be a pair $(H, \alpha)$ where $H$ is an open half-space in $\mathbb{H}^3$ and $\alpha : H \to X$ satisfies

$$\phi \circ \alpha = \mathrm{Id}. \tag{3.2}$$

We call the set $(H_i, \alpha_i)_{i \in I}$ of all half-spaces in $X$ its **half-space decomposition**.

**Lemma 3.2.1**

Let $X$ be a developable hyperbolic end with height function $h$. For all $x \in X$, there exists a unique half-space $(H, \alpha)$ in $X$ such that $x \in \partial\alpha(H)$ and $\nabla h(x)$ is the inward-pointing normal to $\partial\alpha(H)$ at this point.

**Proof:** Let $\phi$ denote a developing map of $X$. Let $x$ be a point of $X$. Define the subset $E^+$ of $T_x X$ by $E^+ := \{\xi \mid \langle \xi, \nabla h(x) \rangle > 0\}$. By Theorem 3.1.1, $E^+$ lies within the domain of the exponential map $\mathrm{Exp}_x$ of $X$ at $x$. By Hadamard's theorem, the composition $(\phi \circ \mathrm{Exp}_x)$ restricts to a diffeomorphism from $E^+$ onto its image $H := (\phi \circ \mathrm{Exp}_x)(E^+)$. This image is an open half-space in $\mathbb{H}^3$ and the function $\alpha := \mathrm{Exp}_x \circ (\phi \circ \mathrm{Exp}_x)^{-1}$ is the desired right-inverse of $\phi$. We readily verify that $(H, \alpha)$ is the desired half-space and that it is unique. This completes the proof. $\square$

**Corollary 3.2.2**

*Let $X$ be a developable hyperbolic end. The half-space decomposition of $X$ covers $X$.*

**Proof:** Indeed, for all $x \in X$, upon applying Lemma 3.2.1 to any point lying vertically below $x$, we obtain a half-space in $X$ containing $x$. The result follows. $\square$

We define the **join-relation** $\sim$ of the half-space decomposition such that, for all $i, j \in I$,

$$i \sim j \ \Leftrightarrow \ \alpha_i(H_i) \cap \alpha_j(H_j) \neq \emptyset. \tag{3.3}$$

This relation is trivially reflexive and symmetric, but not transitive. Composing with $\phi$, we obtain

$$i \sim j \ \Rightarrow H_i \cap H_j \neq \emptyset, \tag{3.4}$$

and

$$i \sim j, \ j \sim k, \ H_i \cap H_j \cap H_k \neq \emptyset \ \Rightarrow \ i \sim k. \tag{3.5}$$

As in Section 2.2, we call the pair $((H_i)_{i \in I}, \sim)$ the **combinatorial data** of $X$. By Theorem 2.2.1 and the subsequent remark, this data is sufficient to recover $X$ up to isometry.

As a first application of the half-space decomposition, we obtain an elementary formula for the height function. Indeed, for all $x \in X$, let $I(x)$ denote the set of indices $i \in I$ such that $x \in \alpha_i(H_i)$.

**Lemma 3.2.3**

*Let $X$ be a developable hyperbolic end with height function $h$, developing map $\phi$, and half-space decomposition $(H_i, \alpha_i)_{i \in I}$. For all $x \in X$,*

$$h(x) = \underset{i \in I(x)}{\mathrm{Sup}} \ d(\phi(x), \partial H_i). \tag{3.6}$$

*In particular, the height function of $X$ is unique.*

**Proof:** Choose $x \in X$. Since the integral curves of the gradient of $h$ are unit speed geodesics,

$$h(x) \geq \underset{i \in I(x)}{\mathrm{Sup}} \ d(\phi(x), \partial H_i).$$

Conversely, by completeness, there exists an integral curve $\gamma :] - h(x), \infty[ \to X$ of $\nabla h$ such that $\gamma(0) = x$. By Lemma 3.2.1, for all $\epsilon > 0$, there exists $k \in I(x)$ such that $\gamma(\epsilon - h(x)) \in \partial \alpha_k(H_k)$ and $\dot{\gamma}(\epsilon - h(x))$ is the inward-pointing normal to $\partial \alpha_k(H_k)$ at this point. In particular,

$$\underset{i \in I(x)}{\mathrm{Sup}} \ d(\phi(x), \partial H_i) \geq d(\phi(x), \partial H_k) = h(x) - \epsilon.$$

Since $\epsilon > 0$ is arbitrary, the result follows. $\square$

**Corollary 3.2.4, Monotonicity**

*Let $X$ and $X'$ be developable hyperbolic ends with respective height functions $h$ and $h'$. If $\psi : X \to X'$ is a morphism, then*

$$h \leq h' \circ \psi. \tag{3.7}$$

**Proof:** Indeed, $\psi$ sends the half-space decomposition of $X$ into the half-space decomposition of $X'$. $\square$

More generally, we obtain the following structure result for morphisms of hyperbolic ends.

**Lemma 3.2.5**

*Let $X$ and $X'$ be developable hyperbolic ends with respective height functions $h$ and $h'$. If $\psi : X \to X'$ is a morphism then, for all $x \in X$,*

$$\langle \nabla h(x), \nabla (h' \circ \psi)(x) \rangle > 0. \tag{3.8}$$

**Proof:** Let $\phi : X \to \mathbb{H}^3$ and $\phi' : X' \to \mathbb{H}^3$ be developing maps such that $\phi = \phi' \circ \psi$. Let $x$ be a point of $X$. By Lemma 3.2.1, there exists a unique half-space $(H, \alpha)$ in $X$ such that $x \in \partial \alpha(H)$ and $\nabla h$ is the inward-pointing normal to $\partial \alpha(H)$ at this point. Observe furthermore that the closure of $\alpha(H)$ is complete in $X$. Its image $(H, \psi \circ \alpha)$ is a half-space in $X'$ such that the closure of $(\psi \circ \alpha)(H)$ is complete in $X'$. Denote $Y' := \partial (\psi \circ \alpha)(H)$ and let $N' : Y' \to TX'$ denote its inward-pointing unit normal vector field. At every point of $Y'$, $\langle N', \nabla h' \rangle > 0$, for otherwise $h'$ would vanish at some point $x'$ of the closure of $(\psi \circ \alpha)(H)$, which is absurd. The result now follows upon pulling back this inequality through $\psi$ and evaluating at $x$. $\square$

**3.3 - Geodesic arcs and convexity.** Geodesic arcs in the half-space decomposition are defined in a similar manner as in the Möbius case. We first consider open half-spaces $H_0$ and $H_1$ in $\mathbb{H}^3$. We say that $H_0$ and $H_1$ **overlap** whenever their boundaries meet along a geodesic. Observe that this holds if and only if their intersection is non-trivial, the intersection of their complements is non-trivial and one is not contained within the other. When $H_0$ and $H_1$ overlap, we define the **geodesic arc** between them to be the set of all open half-spaces $H$ in $\mathbb{H}^3$ such that

$$H_0 \cap H_1 \subseteq H \subseteq H_0 \cup H_1. \tag{3.9}$$

This definition extends to half-spaces in developable hyperbolic ends as follows. Let $X$ be a developable hyperbolic end with developing map $\phi$. We say that two distinct open half-spaces $(H_0, \alpha_0)$ and $(H_1, \alpha_1)$ in $X$ **overlap** whenever the sets $\alpha_0(H_0)$ and $\alpha_1(H_1)$ have non-trivial intersection and neither is contained within the other. Upon composing with $\phi$, it follows that $H_0$ and $H_1$ likewise have non-trivial intersection, and neither is contained within the other. Furthermore, their complements also have non-trivial intersection, for otherwise $X$ would be isometric to $\mathbb{H}^3$, contradicting the existence of a height function. $H_0$ and $H_1$ consequently overlap. Using a connectedness argument, we show that

$$\alpha_0|_{H_0 \cap H_1} = \alpha_1|_{H_0 \cap H_1}, \tag{3.10}$$

so that these functions join to define a function $\alpha_{01} : H_0 \cup H_1 \to X$ such that

$$\phi \circ \alpha_{01} = \mathrm{Id}. \tag{3.11}$$

In particular, for any other open half-space $H$ along the geodesic arc from $H_0$ to $H_1$, $(H, \alpha_{01})$ is also a half-space in $X$. We thus define the **geodesic arc** from $(H_0, \alpha_0)$ to $(H_1, \alpha_1)$ to be the set of all half-spaces in $X$ of this form. We say that a subset $(H_i, \alpha_i)_{i \in J}$ of the half-space decomposition of $X$ is **convex** whenever it contains the geodesic arc between any two of its overlapping elements.

Using this concept of convexity, we obtain deeper information about the structure of the height function. Indeed, let $X$ be a simply connected hyperbolic end with height function $h$, developing map $\phi$ and half-space decomposition $(H_i, \alpha_i)_{i \in I}$. For all $x \in X$, let $I(x)$ be as in the preceding section and observe now that this set is convex. For all $i \in I(x)$, let $r_{x,i}$ denote the supremal radius of open geodesic balls in $H_i$ centred on $\phi(x)$. By Lemma 3.2.3, the height function $h$ of $X$ satisfies

$$h(x) := \sup_{i \in I(x)} r_{x,i}. \tag{3.12}$$

For all $x \in X$, define the partial order $\geq_x$ over $I$ such that, for all $i, j \in I$,

$$i \geq_x j \ \Leftrightarrow i, j \in I(x) \text{ and } r_{x,i} \geq r_{x,j}. \tag{3.13}$$

Define also

$$\hat{\Omega}_x := \bigcup_{i \in I(x)} H_i. \tag{3.14}$$

By a connectedness argument, for all $i, j \in I(x)$,

$$\alpha_i|_{H_i \cap H_j} = \alpha_j|_{H_i \cap H_j}, \tag{3.15}$$

so that these functions join to define a smooth function $\alpha_x : \hat{\Omega}_x \to X$ such that

$$\phi \circ \alpha_x = \text{Id}. \tag{3.16}$$

We call $(\hat{\Omega}_x, \alpha_x)$ the **localisation** of $X$ about $x$. Let $\hat{K}_x$ denote the complement of $\hat{\Omega}_x$ in $\mathbb{H}^3$ and let $h_x : \hat{\Omega}_x \to \mathbb{R}$ denote the distance to $\hat{K}_x$. Since $\hat{K}_x$ is an intersection of closed half-spaces, it is a closed, convex subset of $\mathbb{H}^3$ so that $\hat{\Omega}_x$ is a hyperbolic end with height function $h_x$.

**Lemma and Definition 3.3.1**

*Let $X$ be a developable hyperbolic end with height function $h$ and developing map $\phi$. Let $x$ be a point of $X$, let $(\hat{\Omega}_x, \alpha_x)$ denote the localisation of $X$ at $x$, and let $h_x$ denote its height function. $I(x)$ contains a unique maximal element for $\geq_x$ which realises $h(x)$. Furthermore,*

$$h(x) = (h_x \circ \phi)(x), \tag{3.17}$$

*and for all $y \in \alpha_x(\hat{\Omega}_x)$,*

$$h(y) \geq (h_x \circ \phi)(y). \tag{3.18}$$

*We denote by $\max(x)$ the unique maximal element of $I(x)$.*

**Proof:** Let $x$ be a point of $X$. Since, by (3.12), $(r_{x,i})_{i \in I(x)}$ is bounded above by $h(x)$, $I(x)$ contains a maximal element, and existence follows. Since $I(x)$ is convex and since the restriction of $r_{x,i}$ to every geodesic arc in $I(x)$ is strictly concave, uniqueness follows. Finally, since $\alpha_x$ sends the half-space decomposition of $\hat{\Omega}_x$ to $I(x)$, (3.17) and (3.18) follow, and this completes the proof. $\square$

**Lemma 3.3.2**

*Let $X$ be a developable hyperbolic end with developing map $\phi$. Let $x$ be a point of $X$ and let $(\hat{\Omega}_x, \alpha_x)$ denote the localisation of $X$ about $x$. There exists a neighbourhood $U_x$ of $x$ in $\alpha_x(\hat{\Omega}_x)$ such that, for all $y \in U_x$,*

$$\max(y) \in I(x). \tag{3.19}$$

**Proof:** Let $h$ denote the height function of $X$ and let $(H_i, \alpha_i)_{i \in I}$ denote its half-space decomposition. For $x \in X$, define

$$U_x := \left\{ y \in \alpha_x(\hat{\Omega}_x) \mid d(\phi(y), \phi(x)) < h(x)/2 \right\}.$$

For $y \in U_x$, $h(y) > h(x)/2$. It follows that if $i := \max(y)$, then $H_i$ contains the ball of radius $h(x)/2$ about $\phi(y)$. In particular, $\phi(x)$ is an element of $H_i$, so that $x$ is an element of $\alpha_i(H_i)$, and $i \in I(x)$, as desired. $\square$

**Corollary 3.3.3**

*Let $X$ be a developable hyperbolic end with height function $h$ and developing map $\phi$. Let $x$ be a point of $X$, let $(\hat{\Omega}_x, \alpha_x)$ denote the localisation of $X$ about $x$, and let $h_x$ denote its height function. With $U_x$ as in Lemma 3.3.2, for $y \in U_x$,*

$$h(y) = (h_x \circ \phi)(y). \tag{3.20}$$

**Proof:** Indeed, $\alpha_x$ sends the half-space decomposition of $\hat{\Omega}_x$ to $I(x)$. $\square$

We now use the preceding results to derive more refined analytic properties of the height function. We first require the following definition of PDE theory. Given a smooth manifold $Y$, a point $y \in Y$, a function $f : Y \to \mathbb{R}$ and a symmetric bilinear form $B$ over $T_y Y$, we say that

$$\text{Hess}(f)(x) \geq B \tag{3.21}$$

in the **weak sense** whenever there exists a neighbourhood $\Omega$ of $y$ in $Y$ and a smooth function $g : \Omega \to \mathbb{R}$ such that

(1) $g \leq f$;

(2) $g(y) = f(y)$; and

(3) $\text{Hess}(g)(y) = B$.

We likewise say that $\text{Hess}(f)(y) \leq B$ in the **weak sense** whenever $\text{Hess}(-f)(y) \geq -B$ in the weak sense.

**Lemma 3.3.4**

Let $X$ be a developable hyperbolic end with height function $h$. For all $x \in X$, with respect to the decomposition $T_x X = \mathrm{Ker}(dh(x)) \oplus \langle \nabla h(x) \rangle$,

$$\begin{pmatrix} tanh(h(x))Id & 0 \\ 0 & 0 \end{pmatrix} \leq Hess(h)(x) \leq \begin{pmatrix} coth(h(x))Id & 0 \\ 0 & 0 \end{pmatrix}. \tag{3.22}$$

in the weak sense.

**Proof:** Let $\phi$ denote the developing map of $X$. Let $x$ be a point of $X$. By Corollary 3.3.3, it suffices to prove the result for the localisation $(\hat{\Omega}_x, \alpha_x)$ of $X$ at $x$. Let $y \in \hat{K}_x$ be the closest point to $\phi(x)$. Let $H$ denote the supporting open half-space to $\hat{K}_x$ at $y$ whose boundary is normal to the geodesic joining $y$ to $\phi(x)$. Let $f, g : \mathbb{H}^3 \to \mathbb{R}$ denote respectively the distance to $y$ and the distance to $\partial H$. Trivially, $f(\phi(x)) = g(\phi(x)) = h_x(\phi(x))$ and, over $H$, $g \leq h_x \leq f$. The result now follows upon explicitly determining the hessian operators of $f$ and $g$ at $\phi(x)$. $\square$

**3.4 - Ideal boundaries.** We now study functors which map between the categories of simply connected Möbius surfaces and simply connected hyperbolic ends. We first describe the ideal boundary functor $\partial_\infty$ which associates a simply connected Möbius surface to every simply connected hyperbolic end. For this we require the following finer control of complete geodesic rays in hyperbolic ends.

**Lemma 3.4.1**

Let $X$ be a developable hyperbolic end and let $\gamma : [0, \infty[ \to X$ be a complete geodesic ray. For all $t > 0$, there exists $x \in X_t$ such that

$$\underset{s \to +\infty}{Lim} (\pi_t \circ \gamma)(s) = x. \tag{3.23}$$

In particular, $\gamma$ is asymptotic to the vertical line passing through $x$.

**Proof:** Naturally, we may suppose that $\gamma$ is parametrized by arc-length. Let $h$ denote the height function of $X$. By (3.22), for sufficiently large $t$, the function $f := \langle \dot{\gamma}, \nabla h \circ \gamma \rangle$ satisfies

$$\dot{f}(t) \geq (1 - \epsilon)(1 - f(t)^2)$$

in the weak sense. Solving this inequation, we show that $f$ converges exponentially fast to 1, so that the component of $\dot{\gamma}$ orthogonal to $\nabla h \circ \gamma$ converges exponentially fast to zero. Since $\pi_t$ is 1-Lipschitz, the curve $(\pi_t \circ \gamma)$ thus has finite length, and the result now follows by completeness. $\square$

Let $X$ be a developable hyperbolic end with developing map $\phi$. We define $\partial_\infty X$, the **ideal boundary** of $X$, to be the space of equivalence classes of complete geodesic rays in $X$, where two rays are deemed equivalent whenever they are asymptotic to one another. We furnish this space with the cone topology (see [2]), which is well-defined by Lemma 3.2.1. By Lemma 3.4.1, every complete geodesic ray in $X$ is asymptotic to some vertical line. On the other hand, since $\pi_t$ is 1-Lipschitz for all $t$, no two vertical lines are asymptotic to one another. It follows that $\partial_\infty X$ is diffeomorphic to the leaf space of the vertical line foliation of $X$ which, we recall, is in turn diffeomorphic to every level $X_t$ of $X$. In particular, since $X$ retracts onto $X_t$ for all $t$, it follows that $X$ and $\partial_\infty X$ are homotopy equivalent.

Since the developing map $\phi : X \to \mathbb{H}^3$ sends complete geodesic rays continuously to complete geodesic rays, it defines a continuous function $\partial_\infty \phi : \partial_\infty X \to \partial_\infty \mathbb{H}^3$. By standard properties of convex subsets of hyperbolic space, this function is a local homeomorphism and thus defines a developable Möbius structure over $\partial_\infty X$ which we readily verify is of hyperbolic type.

Finally, let $X'$ be another developable hyperbolic end with developing map $\phi' : X' \to \mathbb{H}^3$ and let $\psi : X \to X'$ be a morphism such that $\phi := \phi' \circ \psi$. Since $\psi$ also maps complete geodesic rays continuously to complete geodesic rays, it defines a morphism $\partial_\infty \psi : \partial_\infty X \to \partial_\infty X'$ such that $\partial_\infty \phi' \circ \partial_\infty \psi = \partial_\infty \phi$. We verify that $\partial_\infty$ respects identity elements and compositions, and thus defines a covariant functor from the category of simply connected hyperbolic ends into the category of simply connected Möbius surfaces.

It is useful to observe how the ideal boundary functor acts on the half-space decomposition of the hyperbolic end.

**Lemma 3.4.2**

*Let $X$ be a developable hyperbolic end with developing map $\phi$, let $(H_i, \alpha_i)_{i \in I}$ denote its half-space decomposition, and let $\sim$ denote its join relation. Then $(\partial_\infty H_i, \partial_\infty \alpha_i)_{i \in I}$ is a subset of the Möbius disk decomposition of $(\partial_\infty X, \partial_\infty \phi)$ which covers $\partial_\infty X$. Furthermore, the restriction to $I$ of the join relation of the Möbius disk decomposition of $\partial_\infty X$ coincides with $\sim$.*

**Remark 3.4.1.** Significantly, however, $(\partial_\infty H_i, \partial_\infty \alpha_i)_{i \in I}$ rarely accounts for the entire Möbius disk decomposition of $\partial_\infty X$. Indeed, this only occurs when $X$ is functionally maximal in the sense of Lemma and Definition 3.5.3, below.

**Proof:** For all $i$, $\partial_\infty H_i$ is a disk in $\hat{\mathbb{C}} = \partial_\infty \mathbb{H}^3$ and, by functoriality, $\partial_\infty \alpha_i$ defines a function from $\partial_\infty H_i$ into $\partial_\infty X$ such that

$$\partial_\infty \phi \circ \partial_\infty \alpha_i = \mathrm{Id}.$$

It follows that $(\partial_\infty H_i, \partial_\infty \alpha_i)_{i \in I}$ is a subset of the Möbius disk decomposition of $\partial_\infty X$. We now show that $(\partial_\infty H_i, \partial_\infty \alpha_i)_{i \in I}$ covers $\partial_\infty X$. Let $\gamma : [0, \infty[ \to X$ be a complete, unit speed geodesic ray. Let $t_0 > 0$ be such that $\dot\gamma(t_0)$ is upward pointing. Let $i$ be the unique element of $I$ such that $\gamma(t_0) \in \partial H_i$ and $\dot\gamma(t_0)$ is the inward pointing unit normal to $\partial H_i$ at this point. The equivalence class of $\gamma$ is then an element of $\partial_\infty \alpha_i(\partial_\infty H_i)$ and since $\gamma$ is arbitrary, it follows that $(\partial_\infty H_i, \partial_\infty \alpha_i)_{i \in I}$ covers $\partial_\infty X$, as desired. Finally, we readily show that, for all $i, j \in I$,

$$i \sim j \;\Leftrightarrow\; \partial_\infty \alpha_i(\partial H_i) \cap \partial_\infty \alpha_j(\partial H_j) \neq \emptyset,$$

so that the restriction to $I$ of the join relation of the Möbius disk decomposition of $\partial_\infty X$ coincides with $\sim$, as desired. $\square$

It is also worth verifying that half-spaces in $X$ are uniquely determined by their ideal boundaries in $\partial_\infty X$.

**Lemma 3.4.3**

*Let $X$ be a developable hyperbolic end. For any two half-spaces $(H_i, \alpha_i)$ and $(H_j, \alpha_j)$ in $X$,*

$$(\partial_\infty H_i, \partial_\infty \alpha_i) = (\partial_\infty H_j, \partial_\infty \alpha_j) \;\Rightarrow\; (H_i, \alpha_i) = (H_j, \alpha_j). \tag{3.24}$$

**Proof:** Since $\partial_\infty H_i = \partial_\infty H_j$, we have $H_i = H_j =: H$. Let $\phi$ denote the developing map of $X$. Denote $U := \alpha_i(H_i) \cap \alpha_j(H_j)$ and $V := \phi(U)$. Observe that, over $V$, $\alpha_i = \phi^{-1} = \alpha_j$. It thus suffices to show that $V = H$. However, since $\alpha_i$ and $\alpha_j$ are local isometries, $\partial V$ is a totally geodesic subset of $\overline{H}$ and, since $\partial_\infty \alpha_i = \partial_\infty \alpha_j$, $\partial V = \partial_\infty H$, so that $V = H$, as desired. $\square$

Finally, the following estimate, though elementary, will play a key role in Chapter 4 in the study of quasicomplete ISC immersions in $\mathbb{H}^3$. Let $X$ be a developable hyperbolic end with developing map $\phi$ and let $(\partial_\infty X, \partial_\infty \phi)$ denotes its ideal boundary. Let $\pi_\infty : X \to \partial_\infty X$ denote the function that sends every point $x \in X$ to the equivalence class of the vertical line on which it lies. We call $\pi_\infty$ the **vertical line projection**.

**Lemma 3.4.4**

*Let $X$ be a developable hyperbolic end with developing map $\phi$, let $(\partial_\infty X, \partial_\infty \phi)$ denote its ideal boundary, let $\omega_\infty$ denote the Kulkarni-Pinkall form of $\partial_\infty X$ and let $\pi_\infty : X \to \partial_\infty X$ denote the vertical line projection. For all $x \in X$,*

$$\phi(x) \in B((\partial_\infty \phi)_*(\omega_\infty \circ \pi_\infty)(x)), \tag{3.25}$$

*where $B$ here denotes the parametrisation of the space of open horoballs in $\mathbb{H}^3$ by $\Lambda^2 \partial_\infty \mathbb{H}^3$ described in Section 1.2.*

**Proof:** Let $h$ denote the height function of $X$, let $x$ be a point of $X$ and denote $x_\infty := \pi_\infty(x)$. Let $y$ be a point of $X$ lying vertically below $x$. In particular, $x_\infty = \pi_\infty(y)$. Let $(H, \alpha)$ be the unique half-space of $X$ such that $y \in \partial \alpha(H)$ and $\nabla h(y)$ is the inward-pointing unit normal to $\partial \alpha(H)$ at this point. Since $\partial_\infty$ is functorial, $(\partial_\infty H, \partial_\infty \alpha)$ is a Möbius disk in $\partial_\infty X$. By definition of the Kulkarni-Pinkall form,

$$(\partial_\infty \phi)_* \omega_\infty(x_\infty) \leq \omega(\partial_\infty H)(\partial_\infty \phi(x_\infty)).$$

Thus, if $B$ is the largest open horoball contained in $H$ with asymptotic centre $\partial_\infty \phi(x_\infty)$, then

$$B \subseteq B((\partial_\infty \phi)_* \omega_\infty(x_\infty)).$$

Since $B$ is an interior tangent to $\partial H$ at $\phi(y)$, it contains $\phi(x)$, and the result follows. $\square$

**3.5 - Extensions of Möbius surfaces.** We now construct the extension functor $\mathcal{H}$, and we will show that it is a right-inverse functor of $\partial_\infty$. Let $S$ be a developable Möbius surface of hyperbolic type with developing map $\phi$. Let $(D_i, \alpha_i)_{i \in I}$ denote its Möbius disk decomposition and let $\sim$ denote its join relation. For all $i$, let $H_i$ denote the open half-space in $\mathbb{H}^3$ with ideal boundary $D_i$. Observe that $((H_i)_{i \in I}, \sim)$ are combinatorial data of some hyperbolic manifold in the sense of Theorem 2.2.1. Let $\mathcal{H}S$, $\mathcal{H}\phi : \mathcal{H}S \to \mathbb{H}^3$ and $(\mathcal{H}\alpha_i)_{i \in I}$ denote respectively the join of $(H_i)_{i \in I}$, its canonical immersion and its canonical parametrisations. In particular, $(H_i, \mathcal{H}\alpha_i)_{i \in I}$ is a half-space decomposition of $\mathcal{H}S$.

In order to show that $\mathcal{H}S$ is a hyperbolic end, it remains only to construct its height function. Bearing in mind Lemma 3.2.3, we proceed as follows. For $x \in \mathcal{H}S$, let $I(x)$ denote the subset of $I$ consisting of those indices for which $x \in \mathcal{H}\alpha_i(H_i)$, and observe that this set is convex. For all $i \in I(x)$, let $r_{x,i}$ denote the supremal radius of open geodesic balls in $H_i$ centred on $\mathcal{H}\phi(x)$. We now define the function $h : \mathcal{H}S \to \mathbb{R}$ by

$$h(x) := \sup_{i \in I(x)} r_{x,i}. \tag{3.26}$$

**Lemma 3.5.1**

*The function $h$ is a height function over $\mathcal{H}S$.*

**Proof:** We first observe that, since $S$ is of hyperbolic type, $I(x)$ contains a maximal element for all $x \in \mathcal{H}S$, and uniqueness of this maximal element follows by the convexity arguments already used earlier in this chapter. The construction and results of Section 3.3 now follow as before. It remains only to verify that $h$ has the required analytic properties. Let $x$ be a point of $\mathcal{H}S$. Let $(\hat{\Omega}_x, \mathcal{H}\alpha_x)$ denote the localisation of $\mathcal{H}S$ at this point and let $h_x$ denote its height function. With $U_x$ as in Lemma 3.3.2, for all $y \in U_x$, $h(y) = (h_x \circ \mathcal{H}\phi)(x)$. It thus follows by standard properties of convex sets in $\mathbb{H}^3$ that $h$ is a locally strictly convex $C_{\mathrm{loc}}^{1,1}$ function whose gradient flow lines are unit speed geodesics. Finally, for all $t > 0$, for all $x \in \mathcal{H}S$ such that $h(x) \geq t$, and for all $\epsilon < t$, the closed ball of radius $(t - \epsilon)$ about $x$ in $\mathcal{H}S$ is complete. It follows that $h^{-1}([t, \infty[)$ is complete for all $t > 0$, and this completes the proof. $\square$

It follows by Lemma 3.5.1 that the operator $\mathcal{H}$ associates a hyperbolic end $\mathcal{H}S$ to every developable Möbius surface $S$. Given two developable Möbius surfaces $S$ and $S'$ and an injective morphism $\phi : S \to S'$, Theorem 2.2.1 yields a canonically defined morphism $\mathcal{H}\phi : \mathcal{H}S \to \mathcal{H}S'$. We verify that $\mathcal{H}$ respects identity elements and compositions and thus defines a covariant functor between these two categories. We call $\mathcal{H}$ the **extension functor**. It is a right inverse of the ideal boundary functor, as the following result shows.

**Lemma 3.5.2**

*Let $S$ be a developable Möbius surface of hyperbolic type with developing map $\phi$, Möbius disk decomposition $(D_i, \alpha_i)_{i \in I}$ and extension $(\mathcal{H}S, \mathcal{H}\phi)$. There exists a unique isomorphism $\psi : S \to \partial_\infty \mathcal{H}S$ such that, for all $i$,*

$$\partial_\infty \mathcal{H}\alpha_i = \psi \circ \alpha_i. \tag{3.27}$$

**Remark 3.5.1.** Naturally, in what follows, rather than mention $\psi$ explicitly, we identify $S$ and $\partial_\infty \mathcal{H}S$.

**Proof:** For all $i$, let $H_i$ denote the open half-space in $\mathbb{H}^3$ with ideal boundary $D_i$. Since $(H_i, \mathcal{H}\alpha_i)_{i \in I}$ is a subset of the half-space decomposition of $\mathcal{H}S$ which covers $\mathcal{H}S$, by Lemma 3.4.2, $(D_i, \partial_\infty \mathcal{H}\alpha_i)_{i \in I}$ is a subset of the Möbius disk decomposition of $\partial_\infty \mathcal{H}S$ which likewise covers $\partial_\infty \mathcal{H}S$. Furthermore, the join relation of this decomposition coincides with that of $(H_i, \mathcal{H}\alpha_i)$, which in turn coincides with that of $(D_i, \alpha_i)_{i \in I}$. It follows by Theorem 2.2.1 that there exists a unique diffeomorphism $\psi : S \to \partial_\infty \mathcal{H}S$ satisfying (3.27), as desired. $\square$

Finally, the height functions of hyperbolic ends obtained by extending Möbius surfaces have more structure than in the general case. Indeed, given a function $f : X \to \mathbb{R}$, a point $x \in X$, a vector $\xi \in T_x(X)$ and a real number $\lambda \in \mathbb{R}$, we say that

$$\mathrm{Hess}(f)(x)(\xi, \xi) \leq \lambda \tag{3.28}$$

in the **weak sense** whenever there exists a geodesic segment $\gamma :] - \epsilon, \epsilon[ \to X$ such that $\gamma(0) = x$, $\dot{\gamma}(x) = \xi$, and

$$\left. \frac{\partial^2}{\partial t^2} f \circ \gamma \right|_{t=0} \leq \lambda. \tag{3.29}$$

in the weak sense of Section 3.3.

**Lemma and Definition 3.5.3**

*Let $S$ be a developable Möbius surface of hyperbolic type, let $\mathcal{H}S$ denote its extension, and let $h$ denote the height function of $\mathcal{H}S$. For all $x \in \mathcal{H}S$, there exists a unit vector $\xi \in T_x\mathcal{H}S$ such that*

$$\langle \xi, \nabla h(x) \rangle = 0, \text{ and}$$
$$Hess(h)(x)(\xi, \xi) \leq tanh(h(x)) \tag{3.30}$$

*in the weak sense. We say that a hyperbolic end $X$ is **functionally maximal** whenever its height function satisfies (3.30).*

**Remark 3.5.2.** We will see in Theorem 3.6.3, below, that a hyperbolic end is an extension of a Möbius surface if and only if it is functionally maximal.

**Proof:** Let $\phi$ denote the developing map of $S$, let $(D_i, \alpha_i)_{i \in I}$ denote its Möbius disk decomposition, and, for all $i$, let $H_i$ denote the open half-space in $\mathbb{H}^3$ with ideal boundary $D_i$. Let $x$ be a point of $\mathcal{H}S$, let $(\hat{\Omega}_x, \mathcal{H}\alpha_x)$ denote the localisation of $\mathcal{H}S$ about $x$ and let $h_x$ denote its height function.

Denote $y := \pi_\infty(x)$. Let $(\Omega_y, \alpha_y)$ denote the localisation of $(S, \phi)$ about $y$ and denote

$$\mathcal{H}\Omega_y := \underset{i \in I(y)}{\cup} H_i.$$

Let $\hat{K}_y$ denote the complement of $\mathcal{H}\Omega_y$ in $\mathbb{H}^3$. Let $h_y : \mathcal{H}\Omega_y \to ]0, \infty[$ denote the distance to $\hat{K}_y$ and observe that $h_y$ is a height function over $\mathcal{H}\Omega_y$ so that $\mathcal{H}\Omega_y$ is a hyperbolic end. Indeed, it is none other than the extension of $\Omega_y$. By functoriality, the extension of $\alpha_y$ is a morphism $\mathcal{H}\alpha_y : \mathcal{H}\Omega_y \to \mathcal{H}S$ such that

$$\mathcal{H}\phi \circ \mathcal{H}\alpha_y = \text{Id}.$$

In particular, $\mathcal{H}\alpha_y$ embeds $\mathcal{H}\Omega_y$ into $\mathcal{H}S$.

For all $i \in I$,

$$x \in \mathcal{H}\alpha_i(H_i) \Rightarrow y \in \alpha_i(D_i),$$

so that every half-space in $(\hat{\Omega}_x, \mathcal{H}\alpha_x)$ is also a half-space in $(\mathcal{H}\Omega_y, \mathcal{H}\alpha_y)$. Consequently,

$$\mathcal{H}\alpha_x(\hat{\Omega}_x) \subseteq \mathcal{H}\alpha_y(\mathcal{H}\Omega_y) \subseteq \mathcal{H}S.$$

It follows by Corollaries 3.2.4 and 3.3.3 that, over $U_x$,

$$h_x \circ \mathcal{H}\phi = h_y \circ \mathcal{H}\phi = h.$$

Now let $z$ denote the closest point in $\hat{K}_y$ to $\mathcal{H}\phi(x)$. Since $\hat{K}_y$ is the convex hull in $\mathbb{H}^3$ of the complement of $\Omega_y$ in $\partial_\infty \mathbb{H}^3$, there exists an open geodesic segment $\gamma :] - \epsilon, \epsilon[ \to \hat{K}_y$ such that $\gamma(0) = z$ (see, for example, [27]). Let $P \subseteq \mathbb{H}^3$ be the totally geodesic plane containing $\gamma$ and $\mathcal{H}\phi(x)$. Let $\xi$ be a unit vector tangent to $P$ at $\mathcal{H}\phi(x)$ and normal to $\nabla h_y$ at this point. We verify that

$$\text{Hess}(h_y)(\mathcal{H}\phi(x))(\xi, \xi) \leq \tanh(x)$$

in the weak sense, and this completes the proof. $\square$

**3.6 - Left inverses and applications.** We study the extent to which $\mathcal{H}$ is also a left inverse of $\partial_\infty$.

**Lemma 3.6.1**

*Let $S$ be a developable Möbius surface with developing map $\phi$. Let $X$ be a developable hyperbolic end with developing map $\psi$. Let $f : \partial_\infty X \to S$ be an injective morphism such that*

$$\phi \circ f = \partial_\infty \psi. \tag{3.31}$$

*There exists a unique injective morphism $\hat{f} : X \to \mathcal{H}S$ such that*

$$\mathcal{H}\phi \circ \hat{f} = \psi, \text{ and}$$
$$\partial_\infty \hat{f} = f. \tag{3.32}$$

**Proof:** Let $(D_i, \alpha_i)_{i \in I}$ denote the Möbius disk decomposition of $S$ with join relation $\sim_\alpha$, let $(H_j, \beta_j)_{j \in J}$ denote the half-space decomposition of $X$ with join relation $\sim_\beta$, and, for all $j \in J$, denote $D_j := \partial_\infty H_j$. By Lemma 3.4.2, $(D_j, \partial_\infty \beta_j)_{j \in J}$ is a subset of the Möbius disk decomposition of $\partial_\infty X$ which covers $\partial_\infty X$. By (3.31), $(D_j, f \circ \partial_\infty \beta_j)_{j \in J}$ is a subset of the Möbius disk decomposition of $S$ which covers $\text{Im}(f)$. We thus identify $J$ with a subset of $I$ in such a manner that, for all $j$,

$$f \circ \partial_\infty \beta_j = \alpha_j.$$

$X$ identifies with the join of $((H_j)_{j \in J}, \sim_\beta)$ whilst the join $Y$ of $((H_j)_{j \in J}, \sim_\alpha)$ identifies with an open subset of $\mathcal{H}S$. However, by injectivity of $f$, the join relations $\sim_\alpha$ and $\sim_\beta$ coincide over $J$, so that, by Theorem 2.2.1, there exists a unique isomorphism $\hat{f} : X \to Y$ such that, for all $j \in J$,

$$\hat{f} \circ \beta_j = \mathcal{H}\alpha_j.$$

Conseqently, for all $j$,

$$\mathcal{H}\phi \circ \hat{f} \circ \beta_j = \mathcal{H}\phi \circ \mathcal{H}\alpha_j = \text{Id} = \psi \circ \beta_j.$$

Likewise, by functoriality,

$$\partial_\infty \hat{f} \circ \partial_\infty \beta_j = \partial_\infty (\hat{f} \circ \beta_j) = \partial_\infty \mathcal{H}\alpha_j = \partial_\infty \mathcal{H}(f \circ \partial_\infty \beta_j).$$

Since $(\beta_j)_{j \in J}$ and $(\partial_\infty \beta_j)_{j \in J}$ cover $X$ and $\partial_\infty X$ respectively, it follows that

$$\mathcal{H}\phi \circ \hat{f} = \psi, \text{ and}$$
$$\partial_\infty \hat{f} = \partial_\infty \mathcal{H}f.$$

Identifying $(\partial_\infty \mathcal{H})f$ with $f$ as in Lemma 3.5.2, we obtain (3.32), and existence follows. To prove uniqueness, let $\hat{f}' : X \to \mathcal{H}S$ be another function satisfying (3.32). Let $j$ be an element of $J$. Since

$$\mathcal{H}\phi \circ \hat{f}' \circ \beta_j = \psi \circ \beta_j = \text{Id} = \mathcal{H}\phi \circ \hat{f} \circ \beta_j,$$

it follows that $(H_j, \hat{f}' \circ \beta_j)$ and $(H_j, \hat{f} \circ \beta_j)$ are half-spaces in $\mathcal{H}S$. Furthermore,

$$\partial_\infty (\hat{f}' \circ \beta_j) = \partial_\infty \hat{f}' \circ \partial_\infty \beta_j = (\partial_\infty \mathcal{H}f) \circ \partial_\infty \beta_j = \partial_\infty \hat{f} \circ \partial_\infty \beta_j = \partial_\infty (\hat{f} \circ \beta_j),$$

so that, by Lemma 3.4.3,

$$\hat{f}' \circ \beta_j = \hat{f} \circ \beta_j.$$

Since $(H_j, \beta_j)_{j \in J}$ covers $X$, it follows that $\hat{f}' = \hat{f}$, and uniqueness follows. $\square$

With Lemma 3.6.1 in mind, we now study the relationship between two hyperbolic ends when one is contained within the other. Thus, let $X$ be a developable hyperbolic end. Let $\mathcal{V}$ denote its vertical line foliation whose leaf space we recall is naturally homeomorphic to $\partial_\infty X$. Let $S$ be a $C^1$ embedded surface in $X$. We say that $S$ is a **graph** over an open subset $\Omega$ of $\partial_\infty X$ whenever it is transverse to $\mathcal{V}$ and the vertical line projection restricts to a homeomorphism from $\partial S$ onto $\Omega$.

**Lemma 3.6.2**

*Let $X$ and $X'$ be developable hyperbolic ends. If $\psi : X \to X'$ is an injective morphism, then the image $\psi(X_t)$ of every level of $X$ is a graph over $\partial_\infty \psi(\partial_\infty X)$ in $X'$.*

**Proof:** Indeed, choose $t > 0$. Let $\pi'_\infty : X' \to \partial_\infty X'$ denote the vertical line projection of $X'$. By Lemma 3.2.5, $Y_t := \psi(X_t)$ is everywhere transverse to the vertical foliation of $X'$ so that the restriction of $\pi'_\infty$ to this surface is everywhere a local homeomorphism. By Theorem 3.1.1, any vertical line which enters $\psi(h^{-1}[t, \infty[)$ remains within this set, so that no vertical line of $X'$ can cross $Y_t$ more than once. It follows that the restriction of $\pi'_\infty$ to this surface is injective.

It only remains to prove surjectivity. By connectedness, it suffices to show that $\pi'_\infty(Y_t)$ is a closed subset of $\partial_\infty \psi(\partial_\infty X)$. Thus, let $(x'_m)_{m \in \mathbb{N}}$ be a sequence of points of $\pi'_\infty(Y_t)$ converging to the limit $x'_\infty \in \partial_\infty \psi(\partial_\infty X)$. For all $m \in \mathbb{N} \cup \{\infty\}$, let $\gamma'_m :]0, \infty[\to X'$ be the height parametrisation of the vertical line of $X'$ terminating at $x'_m$ and, for all finite $m$, let $T_m > 0$ be such that $\gamma'_m(T_m) \in Y_t$. Since $x'_\infty$ is an element of $\partial_\infty \psi(\partial_\infty X)$, there exists $T > 0$ such that $\gamma'_\infty(T) \in \psi(X)$. Since $\psi(X)$ is open, we may therefore suppose that $\gamma'_m(T) \in \psi(X)$ for all $m$. In particular, $T_m < T$ for all $m$, and we may therefore suppose that $(T_m)_{m \in \mathbb{N}}$ converges to some value $T_\infty$, say. For all $m \in \mathbb{N} \cup \{\infty\}$, let $\gamma_m : [T_m - T, \infty[\to X$ denote the unit speed parametrisation of the preimage of $\gamma'_m$ under $\psi$, normalised such that $\psi(\gamma_m(0)) = \gamma'_m(T)$. For all finite $m$, denote $y_m := \mu(T_m - T)$ so that $y_m \in X_t$ and $x'_m = \psi(y_m)$. Since the projection along vertical lines in $X$ to $X_t$ is distance decreasing, it follows that

$$\operatorname*{LimSup}_{m,n \to \infty} d(y_m, y_n) \leq 2T,$$

so that, by completeness, there exists a point $y_\infty$, say, of $X_t$ towards which $(y_m)_{m \in \mathbb{N}}$ subconverges. We verify that $\pi'_\infty(\psi(y_\infty)) = x'_\infty$ so that $\pi'_\infty(Y_t)$ is indeed a closed subset of $\partial_\infty \psi(\partial_\infty X)$. Surjectivity follows, and this completes the proof. $\square$

**Theorem 3.6.3**

*Let $S$ be a developable Möbius surface of hyperbolic type with developing map $\phi$. Let $(\mathcal{H}S, \mathcal{H}\phi)$ denote its extension. $\mathcal{H}S$ is the only functionally maximal developable hyperbolic end with ideal boundary $S$.*

**Proof:** Let $X$ be another developable hyperbolic end with height function $h$ and developing map $\psi$ such that $\partial_\infty X = S$. Let $\hat{f} : X \to \mathcal{H}S$ denote the unique injective morphism such that $\mathcal{H}\phi \circ \hat{f} = \psi$ and $\partial_\infty \hat{f} = \mathrm{Id}$. We identify $X$ with its image $\hat{f}(X)$ in $\mathcal{H}S$. Let $\hat{h}$ denote the height function of $\mathcal{H}S$. By Corollary 3.2.4, $h \leq \hat{h}$. We now claim that $h = \hat{h}$. Indeed, suppose the contrary. Choose $x \in \mathcal{H}S$ such that $\hat{h}(x) > h(x)$. By completeness, for sufficiently small $\epsilon > 0$, there exists a geodesic ray $\gamma : [-h(x) - \epsilon, \infty[\to \mathcal{H}S$ such that $\gamma(0) = x$ and $\dot{\gamma}(0) = \nabla h(x)$. Let $(H, \alpha)$ be the unique half-space in $\mathcal{H}S$ such that $\gamma(-h(x) - \epsilon) \in \partial \alpha(H)$ and $\dot{\gamma}(-h(x) - \epsilon)$ is the inward-pointing unit normal to $\partial \alpha(h)$ at this point, let $h_\epsilon : \alpha(H) \to ]0, \infty[$ denote distance to $\partial \alpha(H)$, and let $Y_\epsilon$ denote the level set of this function at height $h(x) + \epsilon$. For sufficiently small $\epsilon$, $Y_\epsilon$ is wholly contained in $X$ and $h$ restricts to a proper function over this set. Let $y \in Y_\epsilon$ be the point at which $h$ is minimised. Since $x \in Y_\epsilon$, $h(y) \leq h(x)$. However, at this point, with respect to the decomposition $T_y X = \mathrm{Ker}(dh(y)) \oplus \langle \nabla h(y) \rangle$,

$$\mathrm{Hess}(h)(y) \geq \mathrm{Hess}(h_\epsilon)(y) = \begin{pmatrix} \tanh(h(x) + \epsilon) & 0 \\ 0 & 0 \end{pmatrix},$$

which contradicts (3.30). It follows that $h = \hat{h}$ as asserted. Finally, since every level of $X$ is a graph over $S$, it follows that $X = \mathcal{H}S$, as desired. $\square$

We conclude this chapter by addressing the case of non-developable Möbius surfaces and proving Theorem 1.1.4. Observe first that the definition of the ideal boundary function $\partial_\infty$ given in Section 3.4 readily extends to the non-developable case. We now examine the extension functor. Let $S$ be a Möbius surface with fundamental group $\Pi$. Let $\tilde{S}$ denote its universal cover, let $\phi$ be a developing map of $\tilde{S}$ and let $\theta$ denote its holonomy. Let $\mathrm{Deck} : \Pi \to \mathrm{Isom}(\tilde{S})$ denote the action of $\Pi$ on $S$ by deck transformations. By definition, for all $\gamma \in \Pi$,

$$\theta(\gamma) \circ \phi = \phi \circ \mathrm{Deck}(\gamma). \tag{3.33}$$

By Lemma 3.6.1, Deck extends to a unique homeomorphism $\mathcal{H}\mathrm{Deck} : \Pi \to \mathrm{Isom}(\mathcal{H}\tilde{S})$ such that

$$\theta(\gamma) \circ \mathcal{H}\phi = \mathcal{H}\phi \circ \mathcal{H}\mathrm{Deck}(\gamma), \text{ and}$$
$$\partial_\infty(\mathcal{H}\phi \circ \mathcal{H}\gamma) = \phi \circ \mathrm{Deck}(\gamma). \tag{3.34}$$

In addition, for all $\gamma$, $\mathcal{H}\mathrm{Deck}(\gamma)$ preserves every level of $\mathcal{H}\tilde{S}$, and its action on each level is conjugate to its action on $S$. It follows that $\mathcal{H}\mathrm{Deck}$ acts discretely on $\mathcal{H}\tilde{S}$, and we define

$$\mathcal{H}S := \mathcal{H}\tilde{S}/\mathcal{H}\mathrm{Deck}(\Pi). \tag{3.35}$$

We verify that $\mathcal{H}S$ is a hyperbolic end with ideal boundary canonically isomorphic to $S$ and, in the case where $S$ is developable, this hyperbolic end is canonically isomorphic to the extension of $S$ constructed above. This complete the construction of the extension functor in the non-developable case. We now prove Theorem 1.1.4.

**Proof of Theorem 1.1.4:** Suppose that $\mathcal{H}S$ is not maximal. There exists a hyperbolic end $X$ and an injective morphism $f : \mathcal{H}S \to X$ such that $\partial_\infty f : S \to \partial X$ is an isomorphism. In particular, $\partial_\infty X$, and therefore also $X$, has fundamental group $\Pi$. since $\partial_\infty X$ is diffeomorphic to $S$, $\Pi$ is also the fundamental group of $X$. Lifting to the universal covers, $f$ defines a $\Pi$-equivariant injective morphism $\tilde{f} : \mathcal{H}\tilde{S} \to \tilde{X}$ such that $\partial_\infty \tilde{f}$ is an isomorphism. Now let $\phi$ and $\psi$ be respectively developing maps of $\tilde{S}$ and $\tilde{X}$ such that

$$\psi \circ \tilde{f} = \mathcal{H}\phi.$$

Since $\partial_\infty \tilde{f}$ is an isomorphism, by Lemma 3.6.1, there exists a unique injective morphism $g : \tilde{X} \to \mathcal{H}\tilde{S}$ such that $\mathcal{H}\phi \circ g = \psi$ and $\partial_\infty g \circ \partial_\infty \tilde{f} = \mathrm{Id}$. In particular

$$\mathcal{H}\phi \circ (g \circ \tilde{f}) = \mathcal{H}\phi, \text{ and}$$
$$\partial_\infty(g \circ \tilde{f}) = \mathrm{Id},$$

so that, by uniqueness,

$$g \circ \tilde{f} = \mathrm{Id}.$$

Since $g$ is injective, it is also a right-inverse of $\tilde{f}$, so that $\tilde{f}$ is an isomorphism, and therefore so too is $f$. This proves maximality. Uniqueness is proven in a similar manner, and this completes the proof. $\square$

## 4 - Infinitesimally strictly convex immersions.

**4.1 - Infinitesimally strictly convex immersions.** We define an **immersed surface** in $\mathbb{H}^3$ to be a pair $(S, e)$ where $S$ is an *oriented* surface and $e : S \to \mathbb{H}^3$ is a smooth immersion. In what follows, we denote the immersed surface sometimes by $S$ and sometimes by $e$, depending on which is more appropriate to the context. The family of immersed surfaces forms a category where the morphisms between two immersed surfaces $(S, e)$ and $(S, e')$ are those functions $\phi : S \to S'$ such that $e = e' \circ \phi$. Naturally, we identify two immersed surfaces which are isomorphic.

Let $(S, e)$ be an immersed surface. In what follows, we will use the terminology of classical surface theory already described in Section 1.2. Recall that $S$ is said to be **infinitesimally strictly convex** (ISC) whenever its second fundamental form is everywhere positive-definite. When this holds, every point $x$ of $S$ has a neighbourhood $\Omega$ over which $e$ takes values on the boundary of some strictly convex subset $X$ and $N_e$ points *outwards* from this set. Recall also that $S$ is said to be **quasicomplete** whenever the metric $\mathrm{I}_e + \mathrm{III}_e$ is complete. We now show that this is a natural requirement for studying ISC immersions in $\mathbb{H}^3$ in terms of hyperbolic ends. Indeed, denote $\mathcal{E}S := S \times ]0, \infty[$ and define the function $\mathcal{E}e : \mathcal{E}S \to \mathbb{H}^3$ by

$$\mathcal{E}e(x, t) = \mathrm{Exp}(tN_e(x)), \tag{4.1}$$

where Exp here denotes the exponential map of $\mathbb{H}^3$. By standard properties of convex surfaces in $\mathbb{H}^3$, $\mathcal{E}e$ is an immersion and we therefore furnish $\mathcal{E}S$ with the unique hyperbolic metric that makes it into a local isometry.

**Lemma and Definition 4.1.1**

Let $(S, e)$ be an ISC immersed surface in $\mathbb{H}^3$. $\mathcal{E}S$ is a hyperbolic end if and only if $S$ is quasicomplete. When this holds, we call $\mathcal{E}S$ the **end** of $S$, its developing map is $\mathcal{E}e$ and its height function is

$$h(x, t) := t. \tag{4.2}$$

**Proof:** It suffices to verify that $h$ defines a height function over $\mathcal{E}S$ if and only if $S$ is quasicomplete. By definition, $h$ is smooth and its gradient flow lines are unit speed geodesics and, by standard properties of convex subsets of hyperbolic space, $h$ is strictly convex. It thus remains only to study completeness. Choose $t > 0$ and let $e_t$ denote the restriction of $\mathcal{E}e$ to $S \times \{t\}$. By classical hyperbolic geometry, the first fundamental form of this immersion is

$$\mathrm{I}_t := \cosh^2(t)\mathrm{I}_e + 2\cosh(t)\sinh(t)\mathrm{II}_e + \sinh^2(t)\mathrm{III}_e,$$

so that, by infinitesimal strict convexity,

$$\sinh^2(t)(\mathrm{I}_e + \mathrm{III}_e) \le \mathrm{I}_t \le 2\cosh^2(t)(\mathrm{I}_e + \mathrm{III}_e).$$

It follows that $\mathrm{I}_e + \mathrm{III}_e$ is complete if and only if $\mathrm{I}_t$ is complete. However, by convexity, $h^{-1}([t, \infty[)$ is complete if and only if $\mathrm{I}_t$ is complete. Since $t > 0$ is arbitrary, it follows that $h^{-1}([t, \infty[)$ is complete for all $t$ if and only if $S$ is quasicomplete, as desired. $\square$

The operation $\mathcal{E}$ trivially sends morphisms between quasicomplete ISC immersed surfaces to morphisms between hyperbolic ends. Since $\mathcal{E}$ respects identity elements and compositions, it therefore defines a covariant functor between these two categories which we call the **end functor**.

There is also a natural way to associate a Möbius surface to every ISC immersed surface. Indeed, let $(S, e)$ be an ISC immersed surface and let $\phi_e$ denote its asymptotic Gauss map. By infinitesimal strict convexity, $\phi_e$ is a local diffeomorphism from $S$ into $\partial_\infty \mathbb{H}^3 = \hat{\mathbb{C}}$ and is thus the developing map of a Möbius structure over $S$. We denote $\mathcal{M}S := S$ and $\mathcal{M}e := \phi_e$, and we verify that $\mathcal{M}$ defines a covariant functor from the category of ISC immersed surfaces into the category of Möbius surfaces. However, this level of precision will be of little use to us and we will not use this terminology in other sections.

We have now reached a pivotal point of our construction. Indeed, we have associated two distinct hyperbolic ends to each quasicomplete ISC immersed surface, namely the end $\mathcal{E}S$ of $S$ constructed above, and the hyperbolic end $\mathcal{HM}S$ obtained by applying the extension functor of Section 3.5 to the Möbius surface $\mathcal{M}S$. Furthermore, by Lemmas 3.6.1 and 3.6.2, $\mathcal{E}S$ naturally embeds into $\mathcal{HM}S$ in such a manner that the levels of $\mathcal{E}S$ are mapped to graphs in $\mathcal{HM}S$ and, since $e$ is smooth and ISC, a small modification of the proofs of these results extends this embedding to the boundary of $\mathcal{E}S$. This yields an embedding $\tilde{e}$ of $S$ into $\mathcal{HM}S$ which factors the immersion $e$ through the developing map $\mathcal{HM}\phi_e$.

This construction allows us to apply the machinery of Sections 2 and 3 to the study of quasicomplete ISC immersions. Given its utility, we first extend it as follows. Let $(S, \phi)$ be a simply connected Möbius surface and let $(\mathcal{H}S, \mathcal{H}\phi)$ denote its extension. Let $\Omega \subseteq S$ be an open subset and let $e : \Omega \to \mathbb{H}^3$ be a quasicomplete ISC immersion whose Gauss map $\phi_e$ is equal to the restriction of $\phi$ to $\Omega$. Let $(\mathcal{E}\Omega, \mathcal{E}e)$ denote the end of $(\Omega, e)$ and let $\psi : \mathcal{E}\Omega \to \mathcal{H}S$ be the unique injective morphism such that

$$\begin{aligned}
\mathcal{H}\phi \circ \psi &= \mathcal{E}e, \text{ and} \\
\partial_\infty \psi &= \mathrm{Id}.
\end{aligned} \tag{4.3}$$

As before, a small modification of the proof of Lemma 3.6.1 shows that $\psi$ extends to a smooth embedding from the whole of $\overline{\mathcal{E}\Omega} = \Omega \times [0, \infty[$ into $\mathcal{H}S$. We define the embedding $\tilde{e} : S \to \mathcal{H}S$ by

$$\tilde{e}(x) := \psi(x, 0), \tag{4.4}$$

and we call it the **canonical lift** of $e$. By (4.3), $\tilde{e}$ factors $e$ through $\mathcal{H}\phi$ in the sense that

$$e = \mathcal{H}\phi \circ \tilde{e}. \tag{4.5}$$

Furthermore, an equally small modification of the proof of Lemma 3.6.2 then shows that the image of $\tilde{e}$ is also a graph over $\Omega$ in $\mathcal{H}S$, and we denote by $\text{Ext}(\tilde{e})$ the open subset of $\mathcal{H}S$ lying above this graph.

Let $\text{U}^+\mathcal{H}S$ denote the bundle of upward-pointing unit vectors over $\mathcal{H}S$. As in Section 1.2, we define the **horizon map** $\text{Hor} : \text{U}^+\mathcal{H}S \to \partial_\infty \mathcal{H}S = S$ such that, for every unit speed geodesic ray $\gamma : [0, \infty[ \to \mathcal{H}S$ with $\dot{\gamma}(0) = \text{U}^+\mathcal{H}S$,

$$\text{Hor}(\dot{\gamma}(0)) := \lim_{t \to +\infty} \gamma(t). \tag{4.6}$$

This function is well-defined by Theorem 3.1.1. Let $N_{\tilde{e}} : S \to \text{U}^+\mathcal{H}S$ denote the positively-oriented unit normal vector field over $\tilde{e}$. We define the **asymptotic Gauss map** of $\tilde{e}$ by

$$\phi_{\tilde{e}} := \text{Hor} \circ N_{\tilde{e}}. \tag{4.7}$$

**Lemma 4.1.2**

*Let $S$ be a developable Möbius surface of hyperbolic type with developing map $\phi$ and let $(\mathcal{H}S, \mathcal{H}\phi)$ denote its extension. Let $\Omega \subseteq S$ be an open subset and let $e : \Omega \to \mathbb{H}^3$ be a quasicomplete ISC immersion whose asymptotic Gauss map $\phi_e$ is equal to the restriction of $\phi$ to $\Omega$. If $\tilde{e} : \Omega \to \mathcal{H}S$ denotes the canonical lift of $\tilde{e}$, then its asymptotic Gauss map $\phi_{\tilde{e}}$ satisfies*

$$\phi_{\tilde{e}} = \text{Id}. \tag{4.8}$$

**Proof:** Indeed, let $(\mathcal{E}\Omega, \mathcal{E}e)$ denote the end of $(\Omega, e)$. Let $\psi : \overline{\mathcal{E}\Omega} \to \mathcal{H}S$ denote the unique injective morphism such that $\mathcal{H}\phi \circ \psi = \mathcal{E}e$ and $\partial_\infty \psi = \text{Id}$. For all $x \in S$,

$$\phi_{\tilde{e}}(x) = \lim_{t \to \infty} \psi(x, t) = \partial_\infty \psi(x) = x,$$

as desired. $\square$

As a byproduct of the preceding construction, we obtain the following estimate for the Kulkarni-Pinkall metric of $S$ in terms of the geometry of $e$.

**Lemma 4.1.3**

*Let $S$ be a developable Möbius surface with developing map $\phi$ and Kulkarni-Pinkall metric $g_\phi$. Let $\Omega \subseteq S$ be an open subset and let $e : \Omega \to \mathbb{H}^3$ be a quasicomplete ISC immersion whose asymptotic Gauss map $\phi_e$ is equal to the restriction of $\phi$ to $\Omega$. If $I_e$, $II_e$ and $III_e$ denote respectively the first, second and third fundamental forms of $e$ then, over $\Omega$,*

$$g_\phi \leq I_e + 2II_e + III_e. \tag{4.9}$$

**Remark 4.1.1.** The right hand side of (4.9) is the none other than the horospherical metric studied by Schlenker in [23].

**Proof:** Let $(\mathcal{H}S, \mathcal{H}\phi)$ denote the extension of $(S, \phi)$. Let $(\mathcal{E}\Omega, \mathcal{E}e)$ denote the end of $(\Omega, e)$ and let $h$ denote its height function. Let $\psi : \overline{\mathcal{E}\Omega} \to \mathcal{H}S$ denote the unique injective morphism such that $\mathcal{H}\phi \circ \psi = \mathcal{E}e$ and $\partial_\infty \psi = \text{Id}$. Let $x$ be a point of $\Omega$. Let $(H, \alpha)$ denote the half space in $\overline{\mathcal{E}\Omega}$ such that $(x, 0) \in \partial \alpha(H)$ and $\partial_t = \nabla h(x)$ is the inward-pointing normal to this surface at this point. Denote $P := \partial H$ and $D := \partial_\infty H$. Let $g'$ denote the hyperbolic metric of $D$ and let $\phi_0 : P \to D$ denote the asymptotic Gauss map of $P$. Observe that $\phi_0$ is an isometry from $P$ into $(D, g')$.

Let $\text{Hor}$ denote the horizon map of $\text{U}\mathbb{H}^3$. Let $N_e : S \to \text{U}\mathbb{H}^3$ denote the positively-oriented unit normal vector field of $S$. Denote $\nu := N_e(x)$ and observe that, by the chain rule,

$$D\phi(x) = D\text{Hor}(\nu) \circ DN_e(x).$$

Recall now that $T_\nu \text{U}\mathbb{H}^3$ decomposes as

$$T_\nu \text{U}\mathbb{H}^3 = H_\nu \oplus V_\nu, \tag{4.10}$$

where $V_\nu$ is the vertical subspace and $H_\nu$ is the horizontal subspace of the Levi-Civita covariant derivative. Recall furthermore that there exists a natural projection $p_\nu : V_\nu \to \langle\nu\rangle^\perp$ and that $D\pi(\nu)$ maps $H_\nu$ isomorphically onto $T_y\mathbb{H}^3$. We henceforth identify vectors in $V_\nu$ and $H_\nu$ with their respective images under $p_\nu$ and $D\pi(\nu)$. With respect to the decomposition (4.10), for all $\xi \in T_x S$,

$$DN_e(x) \cdot \xi = (De(x) \cdot \xi, De(x) \cdot A_e(x) \cdot \xi),$$

where $A_e$ here denotes the shape operator of $e$. Since $P$ is totally geodesic, for all $\xi \in \langle\nu\rangle^\perp$,

$$D\mathrm{Hor}(\nu) \cdot (\xi, 0) = D\phi_0(y) \cdot \xi.$$

Using the fact that Hor restricts to a conformal diffeomorphism from $U_y\mathbb{H}^3$ into $\partial_\infty\mathbb{H}^3$, we likewise show that, for all $\xi \in \langle\nu\rangle^\perp$,

$$D\mathrm{Hor}(\nu) \cdot (0, \xi) = D\phi_0(y) \cdot \xi.$$

Combining the above relations, it follows that, for all $\xi \in T_x S$,

$$D\phi(x) \cdot \xi = D\phi_0(y) \cdot De(x) \cdot (\xi + A_e(x) \cdot \xi).$$

Since $\phi_0$ is an isometry, it follows that

$$(\phi^* g')(x) = \mathrm{I}_e + 2\mathrm{II}_e + \mathrm{III}_e. \tag{4.11}$$

However, by definition of the Kulkarni-Pinkall metric,

$$g(x) \le (\phi^* g')(x),$$

and the result follows. $\square$

The proof of Lemma 4.1.3 also yields the following useful result.

**Lemma 4.1.4**

*Let $(S, e)$ be an immersed surface in $\mathbb{H}^3$ and let $\mathrm{I}_e$, $\mathrm{II}_e$ and $\mathrm{III}_e$ denote respectively its first, second and third fundamental forms. The asymptotic Gauss map $\phi_e$ of $e$ is conformal with respect to the non-negative semi-definite bilinear form $\mathrm{I}_e + 2\mathrm{II}_e + \mathrm{III}_e$.*

**Proof:** Indeed, this follows from (4.11) since, with the notation of the proof of Lemma 4.1.3, $g'$ is a conformal metric over $D$. $\square$

**4.2 - A priori estimates.** We are now ready to derive our main a priori estimates for quasicomplete ISC immersed surfaces in $\mathbb{H}^3$. First, for every open half-space $H$ in $\mathbb{H}^3$, denote

$$H_r := \{x \in H \mid d(x, \partial H) \ge r\}. \tag{4.12}$$

**Lemma 4.2.1**

*Let $S$ be a developable Möbius surface with developing map $\phi$. Let $(\mathcal{H}S, \mathcal{H}\phi)$ denote its extension, let $(D_i, \alpha_i)_{i \in I}$ denote its half-space decomposition and, for all $i$, let $H_i$ denote the open half-space in $\mathbb{H}^3$ with ideal boundary $D_i$. Let $\Omega$ be an open subset of $S$ and let $e : \Omega \to \mathbb{H}^3$ be a quasicomplete ISC immersion whose asymptotic Gauss map $\phi_e$ is equal to the restriction of $\phi$ to $\Omega$. Let $\tilde{e} : \Omega \to \mathcal{H}S$ denote the canonical lift of $e$ and let $\mathrm{Ext}(\tilde{e})$ denote the subset of $\mathcal{H}S$ lying above $\tilde{e}(\Omega)$. For all $r > 0$, if the extrinsic curvature of $e$ satisfies*

$$K_e \le \tanh(r)^2, \tag{4.13}$$

*then, for all $i$ such that $\alpha_i(D_i) \subseteq \Omega$,*

$$\mathcal{H}\alpha_i(H_{i,r}) \subseteq \mathrm{Ext}(\tilde{e}). \tag{4.14}$$

**Proof:** Let $i$ be an element of $I$ such that $D_i \subseteq \Omega$. Let $j \in I$ be another element such that $\overline{D}_j \subseteq D_i$ and $\alpha_j = \alpha_i|_{D_j}$. Since $\tilde{e}(\Omega)$ is a graph over $\Omega$, for sufficiently large $s$,

$$\mathcal{H}\alpha_j(H_{j,s}) \subseteq \mathrm{Ext}(\tilde{e}).$$

Let $s_0$ be the infimal value of $s$ for which this relation holds. Then the surface $\partial\mathcal{H}\alpha_j(H_{j,s_0})$ is an exterior tangent to $\tilde{e}(\Omega)$ at some point. Since $\partial\mathcal{H}\alpha_j(H_{j,s_0})$ has constant extrinsic curvature equal to $\tanh(s_0)^2$, it follows by the geometric maximum principle that $s_0 \le r$. The result follows upon letting $D_j$ tend to $D_i$. $\square$

**Theorem 4.2.2**

*Let $S$ be a developable Möbius surface with developing map $\phi$. Let $(D_i, \alpha_i)_{i \in I}$ denote its Möbius disk decomposition and, for all $i$, let $H_i$ denote the open half-space in $\mathbb{H}^3$ with ideal boundary $D_i$. Let $\Omega$ be an open subset of $S$ and let $e : \Omega \to \mathbb{H}^3$ be a quasicomplete ISC immersion whose asymptotic Gauss map $\phi_e$ is equal to the restriction of $\phi$ to $\Omega$. For all $r > 0$, if the extrinsic curvature of $e$ satisfies*

$$K_e \leq \tanh(r)^2, \tag{4.15}$$

*then, for all $x \in \Omega$ and for all $i \in I$,*

$$x \in \alpha_i(D_i) \subseteq \Omega \ \Rightarrow e(x) \notin H_{i,r}. \tag{4.16}$$

**Proof:** Let $(\mathcal{H}S, \mathcal{H}\phi)$ denote the extension of $S$, let $(\mathcal{E}\Omega, \mathcal{E}e)$ denote the end of $(\Omega, e)$, and let $\psi : \overline{\mathcal{E}\Omega} \to \mathcal{H}S$ denote the unique injective morphism such that $\mathcal{H}\phi \circ \psi = \mathcal{E}e$ and $\partial_\infty \psi = \mathrm{Id}$. Let $x$ be a point of $\Omega$, and let $i$ be an element of $I$ such that $x \in \alpha_i(D_i) \subseteq \Omega$. Define $\gamma(t) := \mathcal{E}e(x,t)$ and $\mu(t) := \psi(x,t)$, so that $\gamma = \mathcal{H}\phi \circ \mu$. Observe that $\mu$ is the geodesic ray in $\mathcal{H}S$ normal to $\tilde{e}(\Omega)$ at $\tilde{e}(x)$. In particular,

$$\lim_{t \to +\infty} \mu(t) = \partial_\infty \psi(x) = x \in \alpha_i(D_i).$$

Since, by Lemma 4.2.1, $\tilde{e}(\Omega)$ lies outside $\mathcal{H}\alpha_i(H_{i,r})$, it follows by the intermediate value theorem that $\mu$ crosses $\partial \mathcal{H}\alpha_i(H_{i,r})$ at some point. Furthermore, by convexity, $\mu$ crosses this surface transversally from the outside to the inside. Composing with $\mathcal{H}\phi$, it follows that $\gamma$ crosses $\partial H_{i,r}$ transversally from the outside to the inside at some point. However, since $\partial H_{i,r}$ is strictly convex, $\gamma$ can meet this surface no more than once, so that $e(x) = \gamma(0)$ lies outside $H_{i,r}$, as desired. $\square$

**Theorem 4.2.3**

*Let $S$ be a developable Möbius surface with developing map $\phi$ and Kulkarni-Pinkall form $\omega_\phi$. Let $\Omega$ be an open subset of $S$ and let $e : \Omega \to \mathbb{H}^3$ be a quasicomplete ISC immersion whose asymptotic Gauss map $\phi_e$ is equal to the restriction of $\phi$ to $\Omega$. For all $x \in S$,*

$$e(x) \in \overline{B}(\phi_*\omega_\phi(x)), \tag{4.17}$$

*where $B$ here denotes the parametrisation of the space of open horoballs in $\mathbb{H}^3$ by $\Lambda^2 \partial_\infty \mathbb{H}^3$ described in Section 1.2.*

**Proof:** Let $(\mathcal{E}\Omega, \mathcal{E}e)$ denote the end of $(\Omega, e)$. Observe that the level foliation of $(\mathcal{E}\Omega, \mathcal{E}e)$ is $(\Omega \times \{t\})_{t>0}$ so that every level of this hyperbolic end as well as its ideal boundary $\partial_\infty \mathcal{E}\Omega$ naturally identifies with $\Omega$. With respect to these identifications, $\partial_\infty \mathcal{E}e = \phi$ and the vertical line projection $\pi_\infty : \mathcal{E}\Omega \to \partial_\infty \mathcal{E}\Omega$ is given by $\pi_\infty(x,t) = x$. Thus, if $\omega$ denotes the Kulkarni-Pinkall form of $(\Omega, \phi)$ then, by Lemma 3.4.4, for all $x \in \Omega$ and for all $t > 0$,

$$\mathcal{E}e(x,t) \in B(\phi_*\omega(x)).$$

However, by Lemma 2.4.1, $\omega \geq \omega_\phi$, so that, for all such $x$ and $t$,

$$\mathcal{E}e(x,t) \in B(\phi_*\omega_\phi(x)),$$

and the result now follows upon letting $t$ tend to zero. $\square$

**4.3 - Cheeger-Gromov convergence.** In order for the text to be as self-contained as possible, we now recall the basic theory of Cheeger-Gromov convergence. A **pointed riemannian manifold** is a triplet $(X, g, x)$, where $X$ is a smooth manifold, $g$ is a riemannian metric and $x$ is a point of $X$. We say that a sequence $(X_m, g_m, x_m)_{m \in \mathbb{N}}$ of complete pointed riemannian manifolds converges to the complete pointed riemannian manifold $(X_\infty, g_\infty, x_\infty)$ in the **Cheeger-Gromov sense** whenever there exists a sequence $(\Phi_m)_{m \in \mathbb{N}}$ of functions such that

(1) for all $m$, $\Phi_m : X_\infty \to X_m$ and $\Phi_\infty(x_\infty) = x_m$; and

for every relatively compact open subset $\Omega$ of $X_\infty$, there exists $M$ such that

(2) for all $m \geq M$, the restriction of $\Phi_m$ to $\Omega$ defines a smooth diffeomorphism onto its image; and

(3) the sequence $((\Phi_m|_\Omega)^* g_m)_{m \geq M}$ converges to $g_\infty|_\Omega$ in the $C^\infty_{\text{loc}}$ sense.

We call $(\Phi_m)_{m \in \mathbb{N}}$ a sequence of **convergence maps** of $(X_m, g_m, x_m)_{m \in \mathbb{N}}$ with respect to $(X_\infty, g_\infty, x_\infty)$.

At first sight, the concept of Cheeger-Gromov convergence can appear rather daunting and, indeed, its correct usage can be sometimes counterintuitive. However, it is reassuring to observe that it defines a Hausdorff topology over the space of isometry equivalence classes of complete pointed riemannian manifolds.[*] Furthermore, although the convergence maps are trivially non-unique, any two sequences $(\Phi_m)_{m \in \mathbb{N}}$ and $(\Phi'_m)_{m \in \mathbb{N}}$ of convergence maps are equivalent in the sense that there exists an isometry $\Psi : X_\infty \to X_\infty$ preserving $x_\infty$ such that, for any two relatively compact open subsets $U \subseteq \overline{U} \subseteq V$ of $X_\infty$, there exists $M$ such that

(1) for all $m \geq M$, the respective restrictions of $\Phi_m$ and $\Phi'_m \circ \Psi$ to $U$ and $V$ define smooth diffeomorphisms onto their images;

(2) for all $m \geq M$, $(\Phi'_m \circ \Psi)(U) \subseteq \Phi_m(V)$; and

(3) the sequence $((\Phi_m|_V)^{-1} \circ \Phi'_m \circ \Psi)_{m \in \mathbb{N}}$ converges in the $C^\infty$ sense to the identity map over $U$.

The concept of Cheeger-Gromov convergence applies to sequences of immersed submanifolds as follows. We say that a sequence $(S_m, x_m, \phi_m)_{m \in \mathbb{N}}$ of complete pointed immersed submanifolds in a complete riemannian manifold $(X, g)$ converges to the complete pointed immersed submanifold $(S_\infty, x_\infty, \phi_\infty)$ in the **Cheeger-Gromov sense** whenever $(S_m, x_m, \phi_m^* g)_{m \in \mathbb{N}}$ converges to $(S_\infty, x_\infty, \phi_\infty^* g)$ in the Cheeger-Gromov sense and, for one, and therefore for any, sequence $(\Phi_m)_{m \in \mathbb{N}}$ of convergence maps, the sequence $(\phi_m \circ \Phi_m)_{m \in \mathbb{N}}$ converges to $\phi_\infty$ in the $C^\infty_{\text{loc}}$ sense.

Cheeger-Gromov convergence of immersed submanifolds can also be described in terms of graphs. Indeed, let $NS_\infty$ denote the normal bundle of $(S_\infty, \phi_\infty)$ in $\phi_\infty^* TX$. Recall that the exponential map of $X$ defines a smooth function $\text{Exp} : NS_\infty \to X$. In particular, given a sufficently small smooth section $f : \Omega \to NS_\infty$ defined over an open subset $\Omega$ of $S_\infty$, the composition $\text{Exp} \circ f$ defines a smooth immersion of $\Omega$ in $X$ which we call the **graph** of $f$. It is straightforward to show that if the sequence $(S_m, x_m, \phi_m)_{m \in \mathbb{N}}$ converges to $(S_\infty, x_\infty, \phi_\infty)$ in the Cheeger-Gromov sense, then there exists a sequence $(x'_m)_{m \in \mathbb{N}}$ of points in $S_\infty$ and sequences of functions $(f_m)_{m \in \mathbb{N}}$ and $(\alpha_m)_{m \in \mathbb{N}}$ such that

(1) $(x'_m)_{m \in \mathbb{N}}$ converges to $x_\infty$;

(2) for all $m$, $f_m$ maps $S_\infty$ into $NS_\infty$, $\alpha_m$ maps $S_\infty$ into $S_m$ and $\alpha_m(x'_m) = x_m$; and

for every relatively compact open subset $\Omega$ of $S_\infty$, there exists $M$ such that

(3) for all $m \geq M$, $f_m$ restricts to a smooth section of $NS_\infty$ over $\Omega$, $\alpha_m$ restricts to a smooth diffeomorphism of $\Omega$ onto its image, and

$$\text{Exp} \circ f_m|_\Omega = \phi_m \circ \alpha|_\Omega; \text{ and}$$

(4) the sequence $(f_m)_{m \geq M}$ tends to zero in the $C^\infty_{\text{loc}}$ sense.

---

[*] Strictly speaking, of course, the family of all complete pointed riemannian manifolds is not a set. However, by Whitney's theorem, it is sufficient to restrict attention to the family of submanifolds of $\mathbb{R}^m$, for some $m$, with smooth complete metrics defined over them, which is a set.

Finally, these definitions are readily extended in a number of ways. For example, in the case of a sequence $(X_m, g_m, x_m)_{m \in \mathbb{N}}$ of pointed riemannian manifolds, the hypothesis of completeness is unnecessary. Instead, it is sufficient to assume that for all $R > 0$, there exists $M$ such that, for all $m \geq M$, the closed ball of radius $R$ about $x_m$ in $(X_m, g_m)$ is compact. Likewise, in the case of immersed submanifolds, the target space can be replaced with a sequence $(X_m, g_m, x_m)_{m \in \mathbb{N}}$ of pointed riemannian manifolds converging in the Cheeger-Gromov sense to some complete pointed riemannian manifold. Furthermore, it is not necessary to suppose that the riemannian manifolds in this sequence are complete, and so on.

**4.4 - Labourie's theorems and their applications.** For $k > 0$, we define a $k$-**surface** to be a quasicomplete ISC immersed surface in $\mathbb{H}^3$ of constant extrinsic curvature equal to $k$. Let $(S, e)$ be a $k$-surface in $\mathbb{H}^3$. Let $N_e$ denote its positively-oriented unit normal vector field. By quasicompleteness, $N_e$ is a complete immersion of $S$ in $U\mathbb{H}^3$. In order to emphasise our interest in this function as an *immersion* rather than as a *vector field*, following Labourie, we denote $\hat{e} := N_e$ and we call $\hat{e}$ the **Gauss lift** of $e$. In [16], Labourie proves the following result (see also [27]).

**Theorem 4.4.1, Labourie's compactness theorem**

*Choose $k > 0$. Let $(S_m, e_m, x_m)_{m \in \mathbb{N}}$ be a sequence of pointed $k$-surfaces in $\mathbb{H}^3$. For all $m$, let $\hat{e}_m : S_m \to U\mathbb{H}^3$ denote the Gauss lift of $e_m$. If $(\hat{e}_m(x_m))_{m \in \mathbb{N}}$ remains within a compact subset of $\mathbb{H}^3$ then there exists a complete, pointed immersed surface $(S_\infty, \hat{e}_\infty, x_\infty)$ in $U\mathbb{H}^3$ towards which $(S_m, \hat{e}_m, x_m)_{m \in \mathbb{N}}$ subconverges in the Cheeger-Gromov sense.*

Significantly, Theorem 4.4.1 does not affirm that the limit is a lift of some $k$-surface. In order to address this, Labourie introduces what he calls curtain surfaces which are defined as follows. Given a complete geodesic $\Gamma$ in $\mathbb{H}^3$, we denote by $N\Gamma \subseteq U\mathbb{H}^3$ the set of unit normal vectors over $\Gamma$. Observe that $N\Gamma$ is an immersed surface conformally equivalent to $\mathbb{S}^1 \times \mathbb{R}$ with respect to the Sasaki metric of $U\mathbb{H}^3$. We define a **curtain surface** to be any immersed surface $(S, \hat{e})$ in $U\mathbb{H}^3$ which is a cover of $N\Gamma$, for some complete geodesic $\Gamma$.

**Theorem 4.4.2, Labourie's dichotomy**

*Choose $k > 0$. Let $(S_\infty, \hat{e}_\infty)$ be a limit of a sequence of lifts of $k$-surfaces, as in Theorem 4.4.1. If $(S_\infty, \hat{e}_\infty)$ is not a curtain surface, then $(S_\infty, \pi \circ \hat{e}_\infty)$ is a $k$-surface.*

The phenomenon described in Theorem 4.4.2 is best illustrated by the case where $k = 1$. Indeed, by a theorem of Volkov-Vladimirova and Sasaki (see [27]), the only 1-surfaces in $\mathbb{H}^3$ are the horospheres and covers of equidistant cylinders around complete geodesics. Let $\Gamma$ be a complete geodesic in $\mathbb{H}^3$. For all $r > 0$, let $C_r$ denote the cylinder of radius $r$ about $\Gamma$. For all $m$, let $e_m : \mathbb{R}^2 \to C_{1/m}$ be a covering map which is isometric with respect to the sum of its first and third fundamental forms. The sequence $(\mathbb{R}^2, e_m, 0)_{m \in \mathbb{N}}$ subconverges in the $C^\infty_{\text{loc}}$ sense to a smooth function $e_\infty : \mathbb{R}^2 \to \Gamma$. Consequently, viewed as a sequence of pointed immersed surfaces, this sequence degenerates. However, the sequence $(\mathbb{R}^2, \hat{e}_m, 0)_{m \in \mathbb{N}}$ of Gauss lifts converges in the $C^\infty_{\text{loc}}$ sense to a cover of $N\Gamma$, that is, a curtain surface. Labourie's dichotomy affirms that, even in the general case, this is the only mode of degeneration that can occur.

We are now ready to prove one of the main results of this paper.

**Theorem 4.4.3, Monotone convergence**

*Let $S$ be a developable Möbius surface of hyperbolic type with universal cover not isomorphic to $(\mathbb{C}, \text{Exp}(z))$. Let $\phi$ be a developing map of $S$. Let $(\Omega_m)_{m \in \mathbb{N}}$ be a nested sequence of open subsets of $S$ which exhausts $S$. For all $m$, let $e_m : \Omega_m \to \mathbb{H}^3$ be a $k$-surface whose Gauss map $\phi_m$ is equal to the restriction of $\phi$ to $\Omega$. There exists a $k$-surface $e_\infty : S \to \mathbb{H}^3$ towards which $(e_m)_{m \in \mathbb{N}}$ subconverges in the $C^\infty_{\text{loc}}$ sense over $S$.*

**Proof:** For all $m$, let $\hat{e}_m$ denote the Gauss lift of $e_m$. Let $\omega_\phi$ denote the Kulkarni-Pinkall form of $S$. Let $x$ be a point of $S$. We claim that there exists a Möbius disk $(D, \alpha)$ in $S$ such that $x \in \alpha(D)$ and $(\hat{e}_m \circ \alpha)_{m \in \mathbb{N}}$ subconverges in the $C^\infty_{\text{loc}}$ sense. Indeed, let $(D', \alpha')$ be a Möbius disk in $S$ such that $x \in \alpha'(D')$ and the closure of $\alpha'(D')$ in $S$ is compact. By Theorems 4.2.3 and 4.2.2, for all sufficiently large $m$,

$$e_m(x) \in K := \overline{B}(\phi_* \omega_\phi(x)) \setminus \mathcal{H}\alpha'(H'_r),$$

where $H'$ here denotes the open half-space in $\mathbb{H}^3$ with ideal boundary $D'$ and $r := \operatorname{arctanh}(\sqrt{k})$. Since $K$ is compact, it follows by Theorem 4.4.1 that there exists a complete pointed immersed surface $(S_\infty, \hat{e}_\infty, x_\infty)$ in $\mathrm{U}\mathbb{H}^3$ towards which $(\Omega_m, \hat{e}_m, x)_{m \in \mathbb{N}}$ subconverges in the Cheeger-Gromov sense. Denote

$$\phi_\infty := \operatorname{Hor} \circ \hat{e}_\infty,$$

where Hor here denotes the horizon map of $\mathrm{U}\mathbb{H}^3$. Since $\phi_\infty$ is a local diffeomorphism from $S$ into $\partial_\infty \mathbb{H}^3$, it defines a developable Möbius structure over $S_\infty$. Let $(\Phi_m)_{m \in \mathbb{N}}$ be a sequence of convergence maps of $(\Omega_m, \hat{e}_m, x)_{m \in \mathbb{N}}$ with respect to $(S_\infty, \hat{e}_\infty, x_\infty)$. Let $(D'', \alpha'')$ be a Möbius disk about $x_\infty$ in $(S_\infty, \phi_\infty)$ such that $\alpha''(D'')$ is relatively compact in $S_\infty$. Let $M$ be such that, for all $m \geq M$, $\Phi_m$ restricts to a smooth diffeomorphism from $\alpha''(D'')$ onto an open subset $U_m$ of $S$. Since $(\phi_m \circ \Phi_m \circ \alpha'')_{m \geq M}$ converges in the $C^\infty_{\mathrm{loc}}$ sense over $D''$ to $(\phi_\infty \circ \alpha'') = \operatorname{Id}$, upon increasing $M$ and reducing $D''$ if necessary, we may thus suppose that, for all $m \geq M$, $(\phi_m \circ \Phi_m \circ \alpha'')$ is a diffeomorphism onto its image $\Omega_m$ whose inverse we denote by $\beta_m$. Let $D$ be another disk with closure contained in $D''$ such that $x_\infty \in \alpha''(D)$. For sufficiently large $m$, $\Omega_m$ contains $D$ and we therefore define $\alpha_m : D \to D''$ by $\alpha_m := (\Phi_m \circ \alpha'' \circ \beta_m)$. For all such $m$, $(D, \alpha_m)$ is in fact a Möbius disk about $x$ in $(S, \phi)$. In particular, $\alpha_m$ is independent of $m$, and we therefore denote $\alpha := \alpha_m$. By construction $(\hat{e}_m \circ \alpha)_{m \geq M}$ subconverges to $(\hat{e}_\infty \circ \alpha'')$, and it follows that $(D, \alpha)$ is the desired Möbius disk.

A diagonal argument now shows that $(\hat{e}_m)_{m \in \mathbb{N}}$ subconverges in the $C^\infty_{\mathrm{loc}}$ sense to a smooth immersion $\hat{e}_\infty : S \to \mathrm{U}\mathbb{H}^3$ satisfying

$$\operatorname{Hor} \circ \hat{e}_\infty = \phi_\infty.$$

We now claim that $\hat{e}_\infty$ is complete. Indeed, for all $m \in \mathbb{N}$, if $\mathrm{I}_m$, $\mathrm{II}_m$ and $\mathrm{III}_m$ denote respectively the first, second and third fundamental forms of $e_m$, then, by Lemma 4.1.3,

$$\mathrm{I}_m + \mathrm{III}_m \geq \frac{1}{2}(\mathrm{I}_m + 2\mathrm{II}_m + \mathrm{III}_m) \geq g_\phi,$$

where $g_\phi$ here denotes the Kulkarni-Pinkall metric of $(S, \phi)$. However, for all $m$, $\mathrm{I}_m + \mathrm{III}_m$ is also the metric of $\hat{e}_m$. It follows upon taking limits that the metric of $\hat{e}_\infty$ is bounded below by $g$ and is therefore complete, as asserted. Finally, by Labourie's dichotomy, either $e_\infty := \pi \circ e_\infty$ is a $k$-surface, or $\hat{e}_\infty$ is a curtain surface. Since the latter can only occur when $(S, \phi)$ is isomorphic to a cover of $(\mathbb{C}^*, z)$, that is, when its universal cover is isomorphic to $(\mathbb{C}, \operatorname{Exp}(z))$, it follows that $e_\infty$ is a $k$-surface, and this completes the proof. $\square$

**4.5 - Uniqueness and Existence.** We are now ready to prove the main results of this paper. Before proceeding, we require the following technical lemma.

**Lemma 4.5.1**

Let $(S, e)$ be an ISC immersion in $\mathbb{H}^3$ and let $\hat{e} : S \to \mathrm{U}\mathbb{H}^3$ denote its Gauss lift. If $e$ has constant extrinsic curvature equal to $k$ then, for all $t > 0$, the immersion $e_t(x) := \operatorname{Exp}(t\hat{e}(x))$ has curvature at every point strictly greater than $k$ and strictly less than $1$.

**Proof:** Indeed, by the tube formula (see [12]), the shape operator $A_t$ of $e_t$ solves

$$\dot{A}_t = \operatorname{Id} - \mathrm{A}_t^2.$$

Thus, denoting $H_t := \operatorname{Tr}(A_t)$ and $K_t := \operatorname{Det}(A_t)$, we have

$$\frac{\partial}{\partial t} K_t = \operatorname{Tr}(A_t^{-1} - A_t) = \frac{1}{K_t}(1 - K_t)H_t.$$

Solving this ordinary differential equation with $K_0 = k < 1$ yields, for all $t > 0$,

$$k < K_t < 1,$$

as desired. $\square$

The following result is proven by Labourie in [17]. Since its proof fits into the framework developed in this paper, we include it for completeness.

**Theorem 4.5.2, Monotonicity**

*Let $S$ be a developable Möbius surface with developing map $\phi$ and let $(\mathcal{H}S, \mathcal{H}\phi)$ denote its extension. For $i \in \{1, 2\}$, let $\Omega_i \subseteq S$ be an open subset of $S$, let $e_i : \Omega_i \to \mathbb{H}^3$ be a $k$-surface whose asymptotic Gauss map $\phi_i$ is equal to the restriction of $\phi$ to $\Omega_i$, and for each $i$, let $\tilde{e}_i : \Omega_i \to \mathcal{H}S$ denote the canonical lift of $e_i$. If $\Omega_1 \subseteq \Omega_2$ then $\mathrm{Ext}(\tilde{e}_1) \subseteq \mathrm{Ext}(\tilde{e}_2)$.*

**Proof:** Suppose the contrary. Let $U$ denote the set of all points in $\Omega_1$ whose image under $\tilde{e}_1$ lies in the complement of the closure of $\mathcal{H}(\tilde{e}_2)$. For each $i$, let $(\mathcal{E}\Omega_i, \mathcal{E}e_i)$ denote the end of $(\Omega_i, e_i)$ and let $\psi : \overline{\mathcal{E}\Omega_i} \to \mathcal{H}S$ denote the unique injective morphism such that $\mathcal{H}\phi \circ \psi_i = \mathcal{E}e_i$ and $\partial_\infty \psi_i = \mathrm{Id}$. Denote $r := \arctan(\sqrt{k})$. We claim that there exists a unique smooth function $f : U \to [0, r]$ and a unique function $\alpha : U \to \Omega_2$ such that $\alpha$ is a diffeomorphism onto its image and, for all $x \in U$,

$$(\tilde{e}_2 \circ \alpha)(x) = \psi(x, f(x)).$$

Indeed, let $x$ be a point of $U$. Let $(H, \alpha)$ denote the unique half space in $\overline{\mathcal{E}\Omega_1}$ such that $(x, 0) \in \partial \alpha(H)$ and let $D$ denote the ideal boundary of $H$. $(H, \psi_1 \circ \alpha)$ is the unique half-space in $\mathcal{H}S$ which is tangent to $\tilde{e}_1(\Omega)$ at $\tilde{e}_1(x)$. Since

$$\partial_\infty(\psi_1 \circ \alpha)(D) = (\partial \psi_1 \circ \partial_\infty \alpha)(D) = \partial_\infty \alpha(D) \subseteq \Omega_1 \subseteq \Omega_2,$$

it follows by Lemma 4.2.1 that

$$(\psi_1 \circ \alpha)(H_r) \subseteq \mathcal{H}(\tilde{e}_2).$$

In particular, for all $t > r$,

$$\psi_1(x, t) \in \mathcal{H}(\tilde{e}_2).$$

We define $f(x)$ to be the infimal value of $t$ such that $\psi_1(x, t) \in \tilde{e}_2(\Omega_2)$. By Theorem 3.1.1, this is the only value of $t$ such that $\psi_1(x, t) \in \tilde{e}_2(\Omega_2)$. Since, by convexity, the geodesic $t \mapsto \psi_1(x, t)$ is transverse to $\tilde{e}_2(\Omega_2)$ at this point, it follows that $f$ is smooth. the existence of $\alpha$ now follows from the fact that $\tilde{e}_2$ is an embedding, and this proves the assertion.

We now apply a maximum principle at infinity to obtain a contradiction. Let $(x_m)_{m \in \mathbb{N}}$ be a sequence in $U$ such that

$$\lim_{m \to \infty} f(x_m) = f_0 := \sup_{x \in U} f(x) \le r,$$

and, for all $m$, denote $y_m := \alpha(x_m)$. Let $p$ be a fixed point of $\mathbb{H}^3$. For all $m$, let $\beta_m$ be an isometry of $\mathbb{H}^3$ such that $\beta_m(e_1(x_m)) = p$ and, for each $i$, denote $e_{i,m} := \alpha_m \circ e_i$ and let $\hat{e}_{i,m}$ denote its Gauss lift. By Theorem 4.4.1, we may suppose that there exist pointed immersions $(S_{1,\infty}, \hat{e}_{1,\infty}, x_\infty)$ and $(S_{2,\infty}, \hat{e}_{2,\infty}, y_\infty)$ towards which $(\Omega_1, \hat{e}_{1,m}, x_m)_{m \in \mathbb{N}}$ and $(\Omega_2, \hat{e}_{2,m}, y_m)_{m \in \mathbb{N}}$ subconverge in the Cheeger-Gromov sense. Observe that there exist neighbourhoods $U$ of $x_\infty$ in $S_{1,\infty}$, $V$ of $y_\infty$ in $S_{2,\infty}$, a smooth diffeomorphism $\alpha : U \to V$ and a smooth function $f : U \to [0, r]$ such that $\alpha(x_\infty) = y_\infty$, $f$ attains its maximum value of $f_0$ at $x_\infty$ and, for all $z \in U$,

$$(\hat{e}_{2,\infty} \circ \alpha)(z) = \mathrm{Exp}(f(z)\hat{e}_{1,\infty}(z)). \tag{4.18}$$

We now show that this is absurd. Define $e_{1,\infty} : S_{2,\infty} \to \mathbb{H}^3$ by $e_{1,\infty}(z) := \mathrm{Exp}(f_0 \hat{e}_{1,\infty})$. We claim that the extrinsic curvature of this immersion is at every point strictly greater than $k$. Indeed, there are two cases to consider. If $\hat{e}_{1,\infty}$ is a curtain surface, then $e'_1$ is a cylinder of radius $f_0$ about a complete geodesic in $\mathbb{H}^3$ and thus has constant extrinsic curvature equal to 1. On the other hand, if $e_{1,\infty}$ is the lift of a $k$-surface, then, by Lemma 4.5.1, $\hat{e}'_1$ also has extrinsic curvature at every point strictly greater than $k$, and the assertion follows. We now examine the function $e_{2,\infty} := \pi \circ \hat{e}_{2,\infty}$. Once again, there are two cases to consider. If $\hat{e}_{2,\infty}$ is a curtain surface, then $e_{2,\infty}(S_{2,\infty})$ is a complete geodesic $\Gamma$ in $\mathbb{H}^3$ which, by (4.18), is an interior tangent to $e_{1,\infty}$ at $e_{1,\infty}(x_\infty)$. By convexity, this is absurd. Otherwise, $e_{2,\infty}$ is a $k$-surface which is an interior tangent to $e_{1,\infty}$ at $e_{1,\infty}(x_\infty)$, which is also absurd by the geometric maximum principle. We thus obtain a contradiction in all cases, and this completes the proof. $\square$

Theorem 4.5.2 is useful for studying the geometry of $k$-surfaces in $\mathbb{H}^3$. In the present paper, when $\Omega_1 = \Omega_2$, it yields uniqueness.

**Theorem 4.5.3**

*For all $k \in ]0,1[$ and for every developable Möbius surface $S$ with developing map $\phi$, there exists at most one $k$-surface $e : S \to \mathbb{H}^3$ such that $\phi_e = \phi$.*

**Proof:** Indeed, let $e, e' : S \to \mathbb{H}^3$ be $k$-surfaces such that $\phi_e = \phi_{e'} = \phi$. Let $(\mathcal{H}S, \mathcal{H}\phi)$ denote the extension of $(S, \phi)$. Let $\tilde{e}, \tilde{e}' : S \to \mathcal{H}S$ denote the respective canonical lifts of $e$ and $e'$. By Theorem 4.5.2, $\tilde{e}(S) = \tilde{e}'(S)$. From this it follows that $\tilde{e} = \tilde{e}'$ and so $e = e'$, as desired. $\square$

The following result is proven by Labourie in [17] (see also [27]).

**Theorem 4.5.4**

*Let $S$ be a developable Möbius surface with developing map $\phi$. Let $\Omega$ be a relatively compact open subset of $S$ with smooth boundary. There exists a $k$-surface $e : \Omega \to \mathbb{H}^3$ such that $\phi_e = \phi|_\Omega$.*

This yields the following existence result.

**Theorem 4.5.5**

*Let $S$ be a developable Möbius surface of hyperbolic type with developing map $\phi$. If the universal cover of $S$ is not isomorphic to $(\mathbb{C}, \mathrm{Exp}(z))$, then for all $k \in ]0,1[$, there exists a unique $k$-surface $e : S \to \mathbb{H}^3$ such that $\phi_e = \phi$.*

**Proof:** Let $(\Omega_m)_{m \in \mathbb{N}}$ be a nested sequence of relatively compact open subsets of $S$ with smooth boundary which exhausts $S$. By Theorem 4.5.4, for all $m$, there exists a $k$-surface $e_m : \Omega_m \to \mathbb{H}^3$ such that $\phi_{e_m} = \phi|_{\Omega_m}$. By Theorem 4.4.3, there exists a $k$-surface $e : S \to \mathbb{H}^3$ towards which $(e_m)_{m \in \mathbb{N}}$ subconverges in the $C^\infty_{\mathrm{loc}}$ sense. Uniqueness follows by Theorem 4.5.3, and this completes the proof. $\square$

## A - A non-complete $k$-surface.

**A.1 - A non-complete $k$-surface.** In this appendix, we describe a non-complete $k$-surface. We refer the reader to [27] for thorough proofs of the statements made in what follows.

$$f(z) := -\mathrm{Exp}(z)\cosh(z). \tag{A.1}$$

This is the Schwarzian derivative of the function

$$\tilde{\phi}(z) := \mathrm{Exp}(\mathrm{Exp}(z)). \tag{A.2}$$

For $k \in ]0,1[$, let $\tilde{e}_k : \mathbb{C} \to \mathbb{H}^3$ denote the unique $k$-surface solving the asymptotic Plateau problem $(\mathbb{C}, \tilde{\phi})$. By uniqueness, for all $k \in \mathbb{Z}$,

$$\tilde{\phi}(z + 2\pi i k) = \tilde{\phi}(z). \tag{A.3}$$

so that $\tilde{e}_k$ descends to a unique $k$-surface $e_k : \mathbb{C}^* \to \mathbb{H}^3$ such that, for all $z \in \mathbb{C}$,

$$e_k(\mathrm{Exp}(z)) = \tilde{e}_k(z). \tag{A.4}$$

This $k$-surface solves the Plateau problem $(\mathbb{C}^*, \mathrm{Exp}(z))$.

We now identify $\hat{\mathbb{C}}$ with $\partial_\infty \mathbb{H}^3$. In [27], we show that $e_k$ has a cusp at 0 whose end point in $\partial_\infty \mathbb{H}^3 = \hat{\mathbb{C}}$ is 1. We now study the asymptotic geometry of $e_k(z)$ as $z$ tends to infinity in $\mathbb{C}^*$. Let $\Gamma$ denote the geodesic in $\mathbb{H}^3$ joining 0 and $\infty$. For all $y \in \mathbb{R}$, denote

$$L_y := \{x + iy \mid x \in \mathbb{R}\}, \tag{A.5}$$

The image of $L_y$ under $e_k$ converges exponentially fast to a constant speed parametrisation of $\Gamma$ as $y$ tends to $\pm\infty$. On the other hand, the image of $L_y$ under $\mathrm{Exp}$ is a complete radial line rotating at constant speed as $y$ varies. Since, by definition, $\mathrm{Exp}$ is the asymptotic Gauss map of $e_k$, we see that $e_k$ wraps around $\Gamma$, ever tighter, infinitely many times as $y$ tends to $\pm\infty$.

Möbius structures, hyperbolic ends and $k$-surfaces in hyperbolic space.

We now use a heuristic argument to show that $e_k$ is not complete. By Theorem 4.4.3, $e_k$ is the limit as $m$ tends to infinity of the solution $e_{m,k}$ of the asymptotic Plateau problem $(\mathbb{C} \setminus 2m\pi i\mathbb{Z}, \mathrm{Exp}(z))$. However, by uniqueness, for all $m$, $e_{m,k}$ is $2m\pi$ periodic in the $y$ direction with fundamental domain

$$\Omega_m := \{x + iy \mid x \in \mathbb{R}, \ y \in ]-\pi m, \pi m[\} \setminus \{0\} \,. \tag{A.6}$$

In particular, for all $m$, the resulting quotient surface has the conformal type of $\hat{\mathbb{C}} \setminus \{0, 1, \infty\}$ and, by the Gauss-Bonnet Theorem, its metric has area $2\pi/(1 - k)$. Since this area is independent of $m$, upon letting $m$ tend to infinity, it is reasonable to expect that the area induced over $\mathbb{C}^*$ by $e_k$ is also equal to $2\pi/(1 - k)$. Thus, since $e_k$ has the topology of a pointed disk, it cannot be complete, for otherwise its area would be infinite. In fact, we expect the metric induced by $e_k$ over $\mathbb{C}^*$ to be, up to rescaling, isometric to the surface

$$S := \{z := x + iy \in \mathbb{C} \mid y > 0, \ d(z, 2m\mathbb{Z}) > 1\} \,/4m\mathbb{Z}, \tag{A.7}$$

whose fundamental domain is shown in Figure A.1.1.

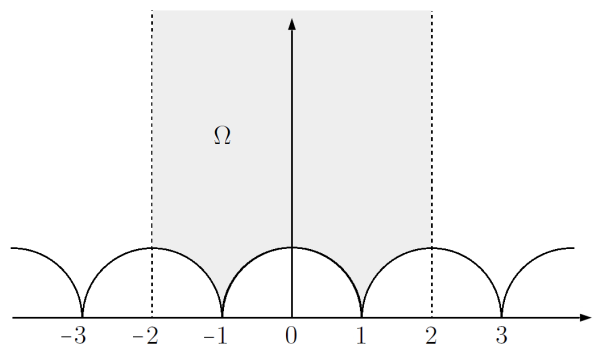

**Figure A.1.1** - The fundamental domain of $e_k$.

# B - Bibliography.

[1] Andersson L., Barbot T., Béguin F., Zeghib A., Cosmological time versus CMC time in spacetimes of constant curvature, *Asian J. Math.*, **16**, no. 2, (2012), 37–87

[2] Ballmann W., Gromov M., Schroeder V., *Manifolds of nonpositive curvature*, Progress in Mathematics, **61**, Birkhäuser-Verlag, Boston, (1985)

[3] Calsamiglia G., Deroin B., Francaviglia S., Branched projective structures with Fuchsian holonomy, *Geom. Topo.*, **18**, (2014), 379–446

[4] Canary R. D., Epstein D., Marden A., *Fundamentals of hyperbolic geometry*, Cambridge University Press, (2006)

[5] Chavel I., *Riemannian geometry, a modern introduction*, Cambridge University Press, (1994)

[6] Dumas D., Complex projective structures, in *Handbook of Teichmller theory. Vol. II*, Eur. Math. Soc., Zrich, 2009, 455–508

[7] Espinar J. M., Gálvez J. A., Mira P., Hypersurfaces in $\mathbb{H}^{n+1}$ and conformally invariant equations: the generalized Christoffel and Nirenberg problems, *J. Eur. Math. Soc.*, **11**, (2009), 903–939

[8] Fillastre F., Seppi A., Spherical, hyperbolic and other projective geometries: convexity, duality, transitions, in *Eighteen Essays in Non-Euclidean Geometry.*, *IRMA Lect. Math. Theor. Phys.*, **29**, Eur. Math. Soc., Zrich, (2019), 321–412

[9] Fillastre F., Smith G., Group actions and scattering problems in Teichmller theory, in *Handbook of group actions IV*, Advanced Lectures in Mathematics, **40**, (2018), 359–417

[10] Gallo D., Kapovich M., Marden A.: The monodromy groups of Schwarzian equations on closed Riemann surfaces, *Ann. Math.*, **151**, no. 2, (2000), 625–704

[11] Goldman W. M., Projective structures with Fuchsian holonomy, *J. Diff. Geom.*, **25**, no. 3, (1987), 297–326

[12] Gromov M., Sign and geometric meaning of curvature, *Rendiconti del Seminario Matematico e Fisico di Milano*, **61**, (1991), 9–123

[13] Kapovich M., *Hyperbolic manifolds and discrete groups*, Birkhaüser-Verlag, (2010)

[14] Kulkarni R. S., Conformal structures and Möbius structures, in *Conformal geometry*, Aspects of Mathematics, MPIM-Bonn, (1988)

[15] Kulkarni R. S., Pinkall U., A canonical metric for Mbius structures and its applications, *Math. Zeit.*, **216**, no. 1, (1994), 89–129

[16] Labourie F., Problèmes de Monge-Ampère, courbes holomorphes et laminations, *Geom. Funct. Anal.*, **7**, no. 3, (1997), 496–534

[17] Labourie F., Un lemme de Morse pour les surfaces convexes, *Invent. Math.*, **141**, (2000), 239–297

[18] Lehto O., *Univalent functions and Teichmüller spaces*, Springer-Verlag, (1990)

[19] Maskit B., *Kleinian groups*, Grundlehren der Mathematischen Wissenschaften, **287**, Springer-Verlag, (1988)

[20] Milnor J., Hyperbolic geometry: the first 150 years, *Bull. Amer. Mat. Soc.*, **6**, no. 1, (1982), 9–24

[21] Santaló L. J., *Geometrías no Euclidianas*, EUDEBA, (1961)

[22] Scannell K. P., Flat conformal structures and the classification of de Sitter manifolds, *Comm. Anal. Geom.*, **7**, no. 2, (1999), 325–345

[23] Schlenker J. M., Hypersurfaces in $\mathbb{H}^n$ and the space of its horospheres, *Geom. Funct. Anal.*, **12**, (2002), 395–435

[24] Smith G., Moduli of Flat Conformal Structures of Hyperbolic Type, *Geom. Dedicata*, **154**, no. 1, (2011), 47–80

[25] Smith G., Hyperbolic Plateau problems, *Geom. Dedicata*, **176**, no. 1, (2014), 31-44

[26] Smith G., *Global Singularity Theory for the Gauss Curvature Equation*, Ensaios Matemáticos, **28**, (2015), 1–114

[27] Smith G., $k$-surfaces in hyperbolic space, *in preparation*

[28] Thurston W. P., Three dimensional manifolds, Kleinian groups and hyperbolic geometry, *Bull. Amer. Mat. Soc.*, **6**, no. 3, (1982), 357–381

[29] Thurston W. P., The geometry and topology of three-manifolds, MSRI, (2002)